\documentclass{amsart}

\usepackage{amssymb,latexsym}
\usepackage{amsthm}
\usepackage{amsmath}
\usepackage[alphabetic]{amsrefs}
\usepackage{ae}
\usepackage{mathrsfs}

\newtheorem{theorem}{Theorem}
\newtheorem{prop}[theorem]{Proposition}
\newtheorem{propdef}[theorem]{Proposition/Definition}
\newtheorem{cor}[theorem]{Corollary}
\newtheorem{lem}[theorem]{Lemma}

\theoremstyle{remark}
\newtheorem{rem}[theorem]{Remark}
\newtheorem{defin}[theorem]{Definition}

\numberwithin{equation}{section}

\renewcommand{\setminus}{\ \rule[2.5pt]{7pt}{.89pt}\ }
\newcommand{\GL}{\operatorname{GL}}

\newcommand{\PGL}{\operatorname{PGL}}
\newcommand{\SL}{\operatorname{SL}}
\newcommand{\Gal}{\operatorname{Gal}}

\newcommand{\gcr}{\ensuremath{G}\text{-cr}}
\newcommand{\up}[2]{{#1^{(#2)}}}

\newcommand{\Hom}{\operatorname{Hom}}
\newcommand{\Der}{\operatorname{Der}}

\newcommand{\tensor}{\otimes}

\newcommand\Z{\mathbf{Z}}    %        integers
\newcommand\F{\mathbf{F}}    %        finite fields
\newcommand\Q{\mathbf{Q}}    %        rational field
\newcommand\G{\mathbf{G}}    %        mult or additive group
\newcommand\Aff{\mathbf{A}}      %        affine space
\newcommand\Lie{\operatorname{Lie}}
\newcommand\Ad{\operatorname{Ad}}
\newcommand\ad{\operatorname{ad}}

\newcommand\lie[1]{\mathfrak{#1}}
\newcommand\glie{\lie{g}}
\newcommand\mlie{\lie{m}}

\newcommand\plie{\lie{p}}

\newcommand{\clie}{\lie{c}}

\newcommand{\zlie}{\lie{z}}

\newcommand{\congruent}{\equiv} 
\newcommand{\iso}{\simeq}

\newcommand{\e}{\varepsilon}

\newcommand{\NN}{\mathcal{N}}
\newcommand{\UU}{\mathcal{U}}

\newcommand{\LL}{\mathcal{L}}

\newcommand{\Int}{\operatorname{Int}}

\newcommand{\sep}{{\operatorname{sep}}}
\newcommand{\Ksep}{K_{\sep}}

\newcommand{\X}{\mathcal{X}}
\newcommand{\T}{\mathcal{T}}

\begin{document}

\author{George J. McNinch}
\address{Department of Mathematics \\ 
         Tufts University \\ 503 Boston Avenue \\ Medford, MA 02155 \\ USA}
\email{george.mcninch@tufts.edu}
\date{21 February, 2005}  

\thanks{Research of the author supported in part by the US National Science Foundation through DMS-0437482}
\title{Optimal $SL(2)$-homomorphisms}


\begin{abstract}
  Let $G$ be a semisimple group over an algebraically closed field of
  \emph{very good} characteristic for $G$.  In the context of
  geometric invariant theory, G. Kempf and -- independently --
  G. Rousseau have associated optimal cocharacters of $G$ to an unstable
  vector in a linear $G$-representation. If the nilpotent element $X
  \in \Lie(G)$ lies in the image of the differential of a homomorphism
  $\SL_2 \to G$, we say that homomorphism is optimal for $X$, or
  simply optimal, provided that its restriction to a suitable torus of
  $\SL_2$ is optimal for $X$ in the sense of geometric invariant theory.
  
  We show here that any two $\SL_2$-homomorphisms which are optimal
  for $X$ are conjugate under the connected centralizer of $X$.  This
  implies, for example, that there is a unique conjugacy class of
  \emph{principal homomorphisms} for $G$. We show that the image of an
  optimal $\SL_2$-homomorphism is a \emph{completely reducible}
  subgroup of $G$; this is a notion defined recently by J-P. Serre.
  Finally, if $G$ is defined over the (arbitrary) subfield $K$ of $k$,
  and if $X \in \Lie(G)(K)$ is a $K$-rational nilpotent element with
  $X^{[p]}=0$, we show that there is an optimal homomorphism
  for $X$ which is defined over $K$.
\end{abstract}

\maketitle
\tableofcontents

\section{Introduction}

Let $G$ be a semisimple group over the algebraically closed field $k$,
and assume that the characteristic of $k$ is \emph{very good} for $G$.
(Actually, we consider in this paper a slightly more general class of
reductive groups; see \S \ref{sec:reductive}, where we also define
\emph{very good} primes).

Premet has recently given a conceptual proof of the Bala-Carter
theorem using ideas of Kempf and of Rousseau from geometric invariant
theory. An element $X \in \glie = \Lie(G)$ is nilpotent just in case
the closure of its adjoint orbit contains 0; such vectors are said to
be unstable.  The Hilbert-Mumford criteria says that an unstable
vector for $G$ is also unstable for certain one-dimensional sub-tori
of $G$. This result has a more precise form due to Kempf and to
Rousseau: there is a class of \emph{optimal} cocharacters of $G$ whose
images exhibit such one dimensional sub-tori. One of the nice features
of these cocharacters is that they each define the same parabolic
subgroup of $G$; for a nilpotent element $X \in \glie$, this
instability parabolic is sometimes called the Jacobson-Morozov
parabolic attached to $X$.

In his proof of the Bala-Carter Theorem in good characteristic,
Pommerening constructed cocharacters associated with the nilpotent
element $X \in \glie$; see \cite{jantzen-nil} for more on this notion,
and see \S \ref{sec:nilpotent} below.  Using the results of Kempf,
Rousseau, and Premet, one finds (cf.  \cite{mcninch-rat}) that the
cocharacters associated with a nilpotent $X \in \glie$ are optimal,
and that any optimal cocharacter $\Psi$ for $X$ such that $X \in
\glie(\Psi;2)$ is associated with $X$ in Pommerening's sense.

In this paper, we show that the notion of optimal cocharacters is
important in the study of subgroups of $G$. We say that a homomorphism
$\phi:\SL_2 \to G$ is optimal provided that the restriction of $\phi$
to the standard maximal torus of $\SL_2$ is a cocharacter associated
to the nilpotent element
\begin{equation*}
  X=  d\phi(
  \begin{pmatrix}
    0 & 1 \\
    0 & 0
  \end{pmatrix}) \in \glie.
\end{equation*}
More precisely, we say that $\phi$ is optimal for $X$.

We prove in this paper that any two optimal homomorphisms for $X$ are
conjugate by $C_G^o(X)$; cf. Theorem \ref{theorem:optimal-conjugate}.
This has an immediate corollary. A principal homomorphism $\phi:\SL_2
\to G$ is one for which the image of $d\phi$ contains a regular
nilpotent element; the conjugacy result just mentioned implies that
there is a unique $G$-conjugacy class of principal homomorphisms.

Generalizing the notion of completely reducible representations, J-P.
Serre has defined the notion of a \gcr~ subgroup $H$ of $G$: $H$ is
\gcr~ if whenever $H$ lies in a parabolic subgroup of $G$, it lies in
a Levi subgroup of that parabolic. We show in Theorem
\ref{theorem:optimal-gcr} that the image of any optimal homomorphism is
\gcr.  In a previous paper \cite{mcninch-sub-principal}, the author
showed the existence of a homomorphism optimal for any $p$-nilpotent $X
\in \glie$; such a homomorphism was essentially obtained (up to
$G$-conjugacy) by base change from a morphism of group schemes defined
over a valuation ring in a number field.  Suppose that $G$ is defined
over the arbitrary subfield $K$ of $k$. If $X$ is a $K$-rational
$p$-nilpotent element, we show in this paper that there is an optimal
homomorphism $\phi$ for $X$ which is defined over $K$; for this we use
the fact, proved in \cite{mcninch-rat}, that some cocharacter
associated with $X$ is defined over $K$.

G. Seitz \cite{seitz} has studied homomorphisms $\phi:\SL_2 \to G$
with the property that all weights of a maximal torus of $\SL_2$ on
$\Lie(G)$ are $\le 2p-2$; he calls the image of such a homomorphism a
good (or restricted) $A_1$-subgroup. We give here a direct proof that
an optimal $\SL_2$-homomorphism is good: we show that the weights of a
cocharacter associated with a $p$-nilpotent element $X \in \glie$ are
all $\le 2p-2$; see Proposition \ref{prop:p-nilpotent-cocharacter}.
It follows from results of Seitz that all good homomorphisms are
optimal -- we do not use this fact in our proofs.

We do use here a result of Seitz (see Proposition \ref{prop:tilting})
to show that $(\Ad \circ \phi,\glie)$ is a tilting module for $\SL_2$
when $\phi$ is the optimal homomorphism obtained previous by the
author \cite{mcninch-sub-principal}; this fact is used to prove a
unicity result Proposition \ref{prop:additive-homom} for certain
homomorphisms $\G_a \to G$ which is crucial to the proof of Theorem
\ref{theorem:optimal-conjugate}; of course, in the end one knows that
$(\Ad \circ \phi,\glie)$ is a tilting module for any optimal $\phi$.

Seitz \emph{loc. cit.} proved a conjugacy result for good
homomorphisms analogous to the result proved here for optimal ones; he
also proved that good homomorphisms are \gcr, so in some sense our
results are not new.  On the other hand, our proofs of conjugacy and
of the \gcr~ property for optimal homomorphisms are free of any case
analysis; we do not appeal to the classification of quasisimple groups
at all. Moreover, we believe that our results on optimal homomorphisms
over ground fields are new and that the ease with which they are obtained
is evidence of the value of our techniques.

As further application of the methods of this paper, we include in \S
\ref{sec:kottwitz} an extension of a result of Kottwitz; we prove that
any nilpotent orbit which is defined over a ground field $K$ contains
a $K$-rational point.

Finally, the appendix contains a note of Jean-Pierre Serre concerning
Springer isomorphisms.

I would like to thank Serre for allowing me to include his note on
Springer isomorphisms as an appendix; I also thank
him for some useful remarks on a preliminary version of this
manuscript. Moreover, I would like to extend thanks to Jens Carsten
Jantzen, and to a referee, for several useful comments on the
manuscript.

\section{Reductive groups}
\label{sec:reductive}

We fix once and for all an algebraically closed field $k$; $K$ will be
an arbitrary subfield of $k$, and $G$ will be a connected, reductive
algebraic group (over $k$) which is defined over the ground field $K$.

If $G$ is quasisimple with root system $R$, the characteristic $p$ of
$k$ is said to be a bad prime for $R$ in the following circumstances:
$p=2$ is bad whenever $R \not = A_r$, $p=3$ is bad if $R =
G_2,F_4,E_r$, and $p=5$ is bad if $R=E_8$.  Otherwise, $p$ is good.
[Here is a more intrinsic definition of good prime: $p$ is good just
in case it divides no coefficient of the highest root in $R$].

If $p$ is good, then $p$ is said to be very good provided that either
$R$ is not of type $A_r$, or that $R=A_r$ and $r \not
\congruent -1 \pmod p$.

If $G$ is reductive, the isogeny theorem \cite{springer-LAG}*{Theorem
9.6.5} yields a -- not necessarily separable -- central isogeny
$\prod_i G_i \times T \to G$ where the $G_i$ are quasisimple and $T$
is a torus. The $G_i$ are uniquely determined by $G$ up to central isogeny,
and $p$ is good (respectively very good) for $G$ if it is good
(respectively very good) for each $G_i$.

The notions of good and very good primes are geometric in the sense
that they depend only on $G$ over $k$. Moreover, they depend only on
the central isogeny class of the derived group $(G,G)$.

We record some facts:
\begin{lem}
  \label{lem:good-char}
  \begin{enumerate}
  \item Let $G$ be a quasisimple group in very good characteristic.
    Then the adjoint representation of $G$ on $\Lie(G)$ is irreducible
    and self-dual. 
  \item Let $M \le G$ be a reductive subgroup containing a maximal
    torus of $G$. If $p$ is good for $G$, then it is good for $M$.
  \end{enumerate}
\end{lem}

\begin{proof}
  For the first assertions of (1), see \cite{Hum95}*{0.13}.  (2) may be found
  for instance in \cite{sommers-mcninch}*{Prop. 16}.
\end{proof}

Consider $K$-groups $H$ which are direct products
\begin{equation*}
  (*)\quad   H = H_1 \times S,
\end{equation*}
where $S$ is a $K$-torus and $H_1$ is a connected, semisimple
$K$-group for which the characteristic is very good. We say that the
reductive $K$-group $G$ is \emph{strongly standard} if there exists a
group $H$ of the form $(*)$ and a separable $K$-isogeny between $G$
and a $K$-Levi subgroup of $H$.  Thus, $G$ is separably isogenous to
$M=C_H(S)$ for some $K$-subtorus $S < H$; note that we do not require
$M$ to be the Levi subgroup of a $K$-rational parabolic subgroup.

We first observe that a strongly standard group $G$ is \emph{standard}
in the sense of \cite{mcninch-rat}; this is contained in the following:

\begin{prop}
  \label{prop:strongly-standard}
  If $G$ is a strongly standard $K$-group, then there is a separable
  $K$-isogeny between $G$ and $\tilde G$ where $\tilde G$ is a reductive $K$-group
  satisfying the ``\emph{standard hypotheses}'' of
  \cite{jantzen-nil}*{\S 2.9}, namely:
  \begin{enumerate}
  \item the derived group of $\tilde G$ is simply connected,
  \item  $p$ is good for $\tilde G$, and
  \item there is a $\tilde G$ invariant nondegenerate bilinear form on
    $\Lie(\tilde G)$.
  \end{enumerate}
\end{prop}

\begin{proof}
  Let $\tilde H = \tilde H_1 \times S$ where $\pi_1:\tilde H_1 \to
  H_1$ is the simply connected cover, and let $\pi = \pi_1 \times
  \operatorname{id}:\tilde H \to H$ be the corresponding isogeny; of
  course, $\tilde H$ and $\pi$ are defined over $H$
  \cite{KMRT}*{Theorem 26.7}.  By assumption, $G=C_H(S)$ for some
  $K$-subtorus $S < H$.  Since $\tilde S = \pi^{-1}(S)^o < \tilde H$
  is again a $K$-torus, its centralizer $\tilde G = C_{\tilde
    H}(\tilde S)$ is a $K$-Levi subgroup of $\tilde H$ and $\pi_{\mid
    \tilde G}:\tilde G \to G$ is an isogeny.  Now, $\Lie(\tilde G)$ is
  the $0$-weight space of $\tilde S$ on $\Lie(\tilde H)$ and $\Lie(G)$
  is the $0$-weight space of $S$ (and $\tilde S$) on $\Lie(H)$.  Since
  $d\pi$ is an $\tilde S$-isomorphism, it restricts to an isomorphism
  $d\pi_{\mid \Lie(\tilde G)} :\Lie(\tilde G) \to \Lie(G)$; in other
  words, $\pi$ is a separable isogeny.
  
  Since $\tilde G$ is a Levi subgroup of $\tilde H$, its derived group
  $\tilde G$ is simply connected, so that (1) holds.  Since $p$ is
  good for $H$, it is also good for $H$ and for the Levi subgroups $G$
  and $\tilde G$; see for instance \cite{sommers-mcninch}*{Prop. 16}.
  Thus (2) holds for $\tilde G$. 
  
  Finally, notice that $\Lie(\tilde H)$ is semisimple as a $\tilde
  H$-module and that $\Lie(H')$ is a self-dual, simple $H'$-module
  whenever $H'$ is quasi-simple in very good characteristic. It
  follows that there is a non-degenerate $\tilde H$-invariant bilinear
  form on $\Lie(\tilde H)$. This restriction of this form to the
  0-weight space for $\tilde S$ is again nondegenerate, and so (3)
  holds.  [Note that the same argument gives non-degenerate invariant
  forms on $\Lie(H)$ and $\Lie(G)$.]
\end{proof}

\begin{rem} Suppose that $V$ is a finite dimensional vector space.  
  Then the group $G=\GL(V)$ is strongly standard.  Indeed, if $\dim V
  \not \congruent 0 \pmod p$, then $G$ is separably isogenous to
  $\SL(V) \times \G_m$, and $p$ is very good for $\SL(V)$.  If $\dim V
  \congruent 0\pmod p$, then $G$ is isomorphic to a Levi subgroup of
  $H = \SL(V \oplus k)$ and $p$ is very good for $H$.
  
  On the other hand, $\SL(V)$ is only strongly standard when $\dim V
  \not \congruent 0 \pmod p$.
\end{rem}

\begin{rem}
  If $G$ is strongly standard, there is always a \emph{symmetric}
  invariant non-degenerate bilinear form on $\Lie(G)$. Indeed, up to
  separable isogeny, $G$ is a Levi subgroup of $T \times H$ where $H$
  is semisimple in very good characteristic. If the result holds for
  $H$, then it holds for $G$; note that any nondegenerate form on
  $\Lie(T)$ is invariant.  Thus we assume that $G$ is semisimple in
  very good characteristic. For such a group, the simply connected
  cover is a separable isogeny so we may also assume $G$ to be simply
  connected.  But then $G$ is a direct product of quasisimple groups,
  hence we may as well suppose that $G$ is quasisimple in very good
  characteristic.  In this case, the adjoint representation is a
  self-dual simple $G$-module. If $p=2$, we are done. Otherwise, one
  can argue as follows: If $G_{/\Q}$ denotes the split group over $\Q$
  with the same root datum as $G$, then the adjoint representation of
  $G_{/\Q}$ is also simple; identifying the weight lattice of a
  maximal torus of $G$ and of $G_{/\Q}$, the adjoint representations
  have the ``same'' highest weight $\lambda$.  Steinberg
  \cite{steinberg-lecture-notes}*{Lemma 79} gives a condition on
  $\lambda$ for the invariant form to be symmetric; since this
  condition is independent of characteristic, and since the Killing
  form is symmetric on $\Lie(G_{/\Q})$, our claim is verified.
%  \footnote{I thank Jens C. Jantzen for pointing out this argument and
%    the reference to Steinberg's lecture notes.}
\end{rem}

\begin{prop}
  \label{prop:separable-orbits}
  If $G$ is strongly standard, then each conjugacy class and each
  adjoint orbit is separable. In particular, if $G$ is
  defined over $K$, and if $g \in G(K)$ and $X \in \glie(K)$, then 
  $C_G(g)$ and $C_G(X)$ are defined over $K$.
\end{prop}

\begin{proof}
  Separability is \cite{springer-steinberg}*{I.5.2 and I.5.6}. The
  fact that the centralizers are defined over $K$ then follows from
  \cite{springer-LAG}*{Prop. 12.1.2}.
\end{proof}

\section{Parabolic subgroups}

In this section, $G$ is an arbitrary reductive group over $k$.
The material we recall here is foundational; the lemmas
from this section will be used mainly for our consideration of
$G$-completely reducible subgroups of a reductive group $G$; cf.
\ref{sub:gcr} below.

If $V$ is an affine variety and $f:\G_m \to V$ is a morphism, we
write $v = \lim_{t\to 0} f(t)$,  and we say that the limit exists, if
$f$ extends to a morphism $\tilde f:k \to V$ with $\tilde f(0) = v$.
If $\gamma$ is any cocharacter of $G$, then
\begin{equation*}
  P_G(\gamma)= P(\gamma) = 
  \{ x \in G \mid \lim_{t \to 0} \gamma(t)x\gamma(t^{-1})\text{\ exists}\}
\end{equation*}
is a parabolic subgroup of $G$ whose Lie algebra is $\plie(\gamma) =
\sum_{i \ge 0} \glie(\gamma;i)$.  Moreover, each parabolic subgroup of
$G$ has the form $P(\gamma)$ for some cocharacter $\gamma$; for all
this cf. \cite{springer-LAG}*{3.2.15 and 8.4.5}.  

We note that $\gamma$ ``exhibits'' a Levi decomposition of
$P=P(\gamma)$. Indeed, $P(\gamma)$ is the semi-direct product
$Z(\gamma) \cdot U(\gamma)$, where $U(\gamma) = \{x \in P \mid
\lim_{t\to 0} \gamma(t)x\gamma(t^{-1}) = 1\}$ is the unipotent radical of
$P(\gamma)$, and the reductive subgroup $Z(\gamma) =
C_G(\gamma(\G_m))$ is a Levi factor in $P(\gamma)$; cf.
\cite{springer-LAG}*{13.4.2}.

\begin{lem}
  \label{lem:parabolic-cocharacter}
  Let $P$ be a parabolic subgroup of $G$, and let $T$ be a maximal
  torus of $P$. Then there is a cocharacter $\gamma \in X_*(T)$ with
  $P=P(\gamma)$.
\end{lem}

\begin{proof}
  Since $P=P(\gamma')$ for some cocharacter $\gamma'$, this follows
  from the conjugacy of maximal tori in $P$.
\end{proof}

For later use, we record:
\begin{lem}
  \label{lem:limit-homomorphism}
  Let $P=P(\gamma)$ be the parabolic subgroup determined by the
  cocharacter $\gamma \in X_*(G)$. Write $L=Z(\gamma)$ for the
  Levi factor of $P$ determined by the choice of $\gamma$. If $\phi:H
  \to P$ is any homomorphism of algebraic groups, the rule
  \begin{equation*}
    \widehat\phi(x) = \lim_{s \to 0} \gamma(s)\phi(x)\gamma(s^{-1})
  \end{equation*}
  determines a homomorphism $\widehat\phi:H \to L$ of algebraic
  groups.  Moreover, the tangent map $d\widehat \phi$ is the composite
  \begin{equation*}
    \Lie(H) \xrightarrow{d\phi} \Lie(P) 
    \xrightarrow{\operatorname{pr}} \Lie(L) = 
    \Lie(P)(\gamma;0)
  \end{equation*}
  where $\operatorname{pr}$ is projection on the 0 weight space.
\end{lem}

\begin{proof}
  It was already observed that $P = L \cdot U$ is a semidirect
  product; the map 
  \begin{equation*}
    x \mapsto \lim_{s \to 0} \gamma(s)x\gamma(s^{-1})
  \end{equation*}
  is the projection of $P$ on $L$ and is thus an algebraic group
  homomorphism $\psi:P \to L$. The tangent map to $\psi$ is evidently
  given by projection onto the 0-weight space for the image of
  $\gamma$, and the lemma follows.
\end{proof}

\begin{rem}
  \label{rem:lim-over-K}
  If the cocharacter $\gamma$ is defined over the ground field $K$,
  then $P = P(\gamma)$ is a $K$-parabolic subgroup, and the Levi
  factor $L= Z(\gamma)$ is defined over $K$. The projection $P \to L$
  given by $x \mapsto \lim_{s \to 0} \gamma(s)x\gamma(s^{-1})$ is of
  course defined over $K$ as well.
\end{rem}

\section{Springer's isomorphisms}

If the characteristic of $k$ is zero, or is ``sufficiently large''
with respect to the group $G$, (some sort of) exponential map defines an
equivariant isomorphism $\exp:\NN \to \UU$ between the nilpotent
variety and the unipotent variety of $G$. Simple examples show the
exponential to be insufficient in general, however, and in 1969, T. A.
Springer \cite{springer-iso} found (the beginnings of) a good
substitute.  See also the outline given in
\cite{springer-steinberg}*{III \S3}. The unipotent variety is known
always to be normal; to make Springer's work complete, one required
also the normality of the nilpotent variety.  Veldkamp obtained that
normality for ``most'' $p$, and Demazure proved it for $G$ satisfying
our hypothesis; cf. \cite{jantzen-nil}*{8.5}.  We summarize
these remarks in the following:

\begin{prop}[Springer]
  \label{prop:springer-isomorphism}
  Let $G$ be a strongly standard $K$-reductive group, where $K$ is any
  subfield of $k$.  There is a $G$-equivariant isomorphism of
  varieties $\Lambda:\UU \to \NN$ which is defined over $K$.
\end{prop}

\begin{proof}[Sketch]
  We just comment briefly on our assumptions on $G$. First, note that
  if $G$ is the direct product of a torus and a semisimple group in very
  good characteristic, there is a separable isogeny $\tilde G \to G$
  where $\tilde G$ is the direct product of a $K$-torus and a simply
  connected semisimple $K$-group in (very) good characteristic.  Moreover,
  the separable isogeny is defined over $K$ and induces equivariant
  $K$-isomorphisms $\tilde \UU \to \UU$ and $\tilde \NN \to \NN$
  (using some hopefully obvious notation); see
  \cite{mcninch-sub-principal}*{Lemma 27}. Now, Springer proved the
  proposition holds for $\tilde G$ -- see the above references-- and
  thus the result for $G$ is true in this case.
  
  Repeating the above argument, we may replace $G$ by a separably
  isogenous group, and thus we suppose that $G = C_H(S)$, where $S$ is
  a $K$-torus in a $K$-group $H$ as in $(*)$ of section \S
  \ref{sec:reductive}; the above remarks show that there is an
  $H$-equivariant isomorphism $\Lambda_H:\UU_H \to \NN_H$ between the
  unipotent and nilpotent varieties for $H$. Since $\UU = (\UU_H)^S$
  and $\NN=(\NN_H)^S$, it is clear that $\Lambda_H\mid_{\UU}$ defines
  the required isomorphism for the varieties associated with $G$.
\end{proof}

\begin{rem}
  \label{rem:springer-iso-restriction}
  Suppose that $\Lambda:\UU \to \NN$ is an equivariant isomorphism
  defined over $K$.  If $P \le G$ is a $K$-parabolic subgroup, Lemma
  \ref{lem:parabolic-cocharacter} makes clear that the restriction
  $\Lambda_{\mid U}:U \to \Lie(U)$ is a $P$-equivariant isomorphism.
  Similarly, if $L \le G$ is a $K$-Levi subgroup, then $\Lambda_{\mid
    \UU_L}:\UU_L \to \NN_L$ is an $L$-equivariant isomorphism.
\end{rem}

The isomorphism $\Lambda$ of the proposition is \emph{quite far} from 
being unique; cf. the appendix of J-P. Serre below.
We summarize the result of that appendix with the following statements,
which we make only in the ``geometric'' setting -- i.e. over $k$ rather than $K$.
\begin{prop}[Serre]
  \label{prop:serre-springer-iso}
  Let $G$ be a strongly standard reductive $k$-group.
  \begin{enumerate}
  \item Fix a regular nilpotent $X \in \glie$. For each regular
    unipotent $v \in C_G(X)$, there is a unique $G$-equivariant
    isomorphism of varieties $\Lambda_v:\UU \to \NN$ with $\Lambda_v(v)
    = X$.
  \item Any two $G$-equivariant isomorphisms $\Lambda,\Lambda':\UU \to
    \NN$  induce the same map on the finite sets of orbits.
  \end{enumerate}
\end{prop}

%\begin{rem}
%  Suppose that $G$ if $K$-quasi-split, i.e. that $G$ contains a Borel
%  subgroup over $K$. Since $p$ is good, $G$ contains a regular
%  nilpotent element $X \in \glie(K)$; cf.  \cite{springer-iso}*{Prop.
%    3.4}. If $\Lambda:\UU \to \NN$ is an isomorphism over $K$, then of
%  course $\Lambda = \Lambda_v$ where $v=\Lambda^{-1}(X)$. Moreover, $v
%  \in G(K)$.  Conversely, for each regular unipotent $v \in
%  C_G(X)(K)$, it is clear that (with the above notation) $\Lambda_v$
%  is an isomorphism defined over $K$.
  
%  When $G$ is not $K$-quasisplit, I don't know if there is a $K$-variety
%  whose rational points parametrize the Springer isomorphisms which
%  are defined over $K$.
%\end{rem}

\section{Frobenius twists and untwists}

Let $K'$ be a perfect field of characteristic $p>0$, and let $K'
\subset K$ be an arbitrary extension of $K'$. We fix an algebraically
closed field $k$ containing $K$.

In this section, algebras are always assumed to be commutative.
Consider a $K'$-algebra $A$. For $r \in \Z$, we may consider the
$K'$-algebra $\up{A}{r}$ which coincides with $A$ as a ring, but where
each $b \in K'$ acts on $\up{A}{r}$ as $b^{p^{-r}}$ does on $A$.  For
an extension field $K$ of $K'$, we write $\up{A}{r}_{/K}$ and $A_{/K}$
for the $K$-algebras obtained by base-change; thus e.g.  $A_{/K} = A
\tensor_{K'} K$.

Let $r \ge 0$ and let $q=p^r$.  There is a $K'$-algebra homomorphism
$F^r:\up{A}{r} \to A$ given by $x \mapsto x^q$. We write $A^q = \{f^q
\mid f \in A\}$; $A^q$ is a $K'$-subalgebra of $A$, and the image of
$F^r$ coincides with $A^q$.

Let $A$ be a $K'$-algebra and an integral domain. We clearly have:
\begin{lem}
  If $r \ge 0$, and $q=p^r$, then $F^r:\up{A}{r} \to A^q$ is an
  isomorphism of $K'$-algebras.
\end{lem}

Write $B=A_{/K}$. Let us notice that $K[B^q] = K[A^q]$.  For $r \ge
0$, consider the algebra homomorphism $F^r_{/K}:\up{A}{r}_{/K} \to
K[A^q] \subset A_{/K}$  given on pure tensors by $f \tensor \alpha
\mapsto f^q \cdot \alpha$ for $f \in \up{A}{r}$ and $\alpha \in K$. We
have more generally
\begin{lem}
  \label{lem:Frob-twist-q-power}
  For $r \ge 0$, $F^r_{/K}:\up{A}{r}_{/K} \to
  K[B^q]$ is an isomorphism, where again $q=p^r$.
\end{lem}

\begin{proof}
  We have observed already that $C = K[B^q]=K[A^q]$ is the $K$-algebra
  generated by $A^q$. According to the previous lemma, the image of
  the restriction of $F^r_{/K}$ to $\up{A}{r} \tensor 1$ is the set of
  $K$-algebra generators $A^q$ of $C$; this implies that $F^r_{/K}$ is
  surjective.
  
  Since $A$ is a domain, the homomorphism $F^r:\up{A}{r} \to A$ is
  injective.  This implies the injectivity of $F^r_{/K}$ 
  since $K$ is flat over $K'$.
\end{proof}

\begin{lem}
  \label{lem:integral}
  Assume that $A$ is \emph{geometrically irreducible}, i.e. that
  $A_{/k}$ is a domain.  Also assume $A$ to be \emph{geometrically
    normal}, i.e.  that $A_{/k}$ is integrally closed in its field of
  fractions $E$.
  %Write $F$ for the field of fractions of $A_{/K}$ and 
  %Write $E$ for the field of fractions of $A_{/k}$. 
  Let $q=p^r$ for $r
  \ge 0$, and let $f \in A_{/K}$.  Then $f \in K[A^q]$ if and only if
  $f \in E^q$.
\end{lem}

\begin{proof}
  We have clearly the implication $\implies$.  Now suppose that $f \in
  E^q$, say $f = g^q$ for $g \in E$. The normality of $A_{/k}$ shows
  then that $g \in A_{/k}$.  We may find $\alpha_1,\dots,\alpha_n \in
  k$ and elements $f_1,\dots,f_n \in A$ such that $g = \sum_{i=1}^n
  \alpha_i f_i$; we may assume as well that $\{f_i \mid 1\le i \le
  n\}$ is a $K'$-linearly independent set.  Since $K'$ is perfect,
  $\{f_i^q \mid 1\le i \le n\}$ is again $K'$-linearly independent.
  Since $f = g^q = \sum_{i=1}^n \alpha_i^q f_i^q \in A_{/K}$, it
  follows that $\alpha_i^q \in K$ for $1 \le i \le n$ and the proof of
  $\Longleftarrow$ is complete.
\end{proof}

\begin{rem}
  It can happen that $A_{/K}$ is a normal domain, but that
    $A_{/k}$ is not normal; cf. \cite{bourbaki-comm}*{exerc.
      V.\S1.23(b)}.
\end{rem}

\begin{lem}
  \label{lem:diff-vanish}
  Let $X$ and $Y$ be irreducible affine $k$-varieties, and let $f:X
  \to Y$ be a dominant morphism. Then the following are equivalent:
  \begin{enumerate}
  \item[(a)] there is a non-empty open subset $W \subset X$ such that
    $df_x \not= 0$ for all $x \in W(k)$.
  \item[(b)] $f^*(k(Y))$ is not contained in $k(X)^p$.
  \end{enumerate}
\end{lem}

\begin{proof}
  For an affine $k$-variety $Z$, let $\Omega_{Z} = \Omega_{k[Z]/k}$ be
  the module of differentials.  The map $f:X \to Y$ determines a map
  $\phi:\Omega_{Y} \to \Omega_{X}$ of $k[Y]$ modules and -- since $f$
  is dominant -- a map $\psi:\Omega_{k(Y)/k} \to \Omega_{k(X)/k}$ of
  $k(Y)$-vector spaces.
  
  It follows from \cite{springer-LAG}*{Theorem 4.3.3} that there are
  non-empty affine open subsets $U$ of $X$ and $V$ of $Y$ such that
  $f$ restricts to a morphism $U \to V$, $\Omega_U$ is a free
  $k[U]$-module of rank $\dim X$, and $\Omega_V$ is a free
  $k[V]$-module of rank $\dim Y$.  Now, $\phi$ restricts to a map
  $\phi_{\mid \Omega_V}:\Omega_V \to \Omega_U$ of $k[V]$-modules, and
  it is clear that $\phi_{\mid \Omega_V} = 0$ if and only if $\psi=0$
  [use that $\Omega_{k(X)/k} = k(X) \tensor_{k[U]} \Omega_U$ together
  with the corresponding statement for $Y$].
  
  Choosing bases of the free modules $\Omega_U$ and $\Omega_V$,
  $\phi_{\mid \Omega_V}$ is given on $\Omega_V$ by a matrix $M$ with
  entries in $k[U]$. For $x \in U(k)$, the map $df_x:T_xU \to
  T_{f(x)}V$ identifies with the map
  \begin{equation*}
    \Hom_{k[U]}(\Omega_U,k_x) \to \Hom_{k[V]}(\Omega_V,k_{f(x)})
  \end{equation*}
  deduced from $\phi_{\mid \Omega_V}$. The open subset of $U$
  defined by the condition $M_x \not = 0$ is non-empty if and only
  $\phi_{\mid \Omega_V} \not = 0$; thus (a) is equivalent to
  the statement $\psi \not = 0.$
  
  Applying \cite{springer-LAG}*{Theorem 4.2.2}, one knows that the
  restriction mapping
  \begin{equation*}
      \Der_k(k(X),k(X)) \to \Der_k(f^*k(Y),k(X))
  \end{equation*}
  is dual to the mapping $\psi:\Omega_{k(Y)/k} \to \Omega_{k(X)/k}$;
  in particular, this restriction is 0 if and only if $\psi=0$.
  
  Now, it is proved for instance in \cite{lang}*{VIII, Prop. 5.4} that
  $z \in k(X)$ is contained in $k(X)^p$ if and only if $D(z)=0$ for
  each $D \in \Der_k(k(X),k(X))$.  The
  assertion (a) $\iff$ (b) follows at once.
\end{proof}

If $X$ is an affine $K'$-variety and $A = K'[X]$, then for $r \in \Z$
we write $\up{X}{r}$ for the $K'$-variety
$\operatorname{Spec}(\up{A}{r})$.  For an arbitrary $K'$-variety $X$,
one defines the $K'$-variety $\up{X}{r}$ by gluing together the
$K'$-varieties $\up{U_i}{r}$ from an affine open covering $\{U_i \mid
1 \le i \le n\}$ of $X$; this construction is independent of the
choice of the covering.

Let $r \ge 0$.  When $X$ is affine, the $r$-th \emph{Frobenius
  morphism} $F^r_X:X \to \up{X}{r}$ is defined to have 
comorphism $F^r:\up{A}{r} \to A$. For an arbitrary $K'$
variety $X$, there is a unique morphism $F^r_X:X \to \up{X}{r}$ whose
restriction to each affine open subset $U$ of $X$ is given by $F^r_U$.

We write $\up{X}{r}_{/K}$ for the base change of the $K'$-variety
$\up{X}{r}$ to $K$.  

\begin{theorem}
  \label{theorem:untwist}
  Let $X$ and $Y$ be geometrically irreducible $K$-varieties.  Assume
  that $X$ is defined over $K'$ and is geometrically normal -- i.e.
  $X_{/k}$ is normal.  Suppose that $f:X \to Y$ is a $K$-morphism
  whose image  contains a positive dimensional sub-variety of
  $Y$.  There is a unique $r \ge 0$ and a unique $K$-morphism
  $g:\up{X}{r}_{/K} \to Y$ such that
  \begin{enumerate}
  \item $f = g \circ F^r_X$, and
  \item there is a non-empty open subset $U$ of $\up{X}{r}$ such that
    $dg_x \not = 0$ for $x \in U(k)$.
  \end{enumerate}
\end{theorem}

\begin{rem}
  \begin{enumerate}
  \item[(a)] Of course, the image of $f$ contains a non-empty open
    subset $U$ of its closure $\overline{f(X)}$
    \cite{springer-LAG}*{Theorem 1.9.5}, so the dimension assumption
    made in the theorem is equivalent to: $U$ has positive
    dimension.
  \item[(b)] The theorem has been known for a long time, but it seems
    to be difficult to give a reference. It was used for instance by
    J-P. Serre in his classification of the inseparable isogenies of
    height 1 of a group variety (and especially of an abelian
    variety), cf. Amer. J. Math. 80 (1958), pp.715-739, sect. 2.
  \end{enumerate}
\end{rem}

\begin{proof}
  Notice that if the theorem is proved when $X$ and $Y$ are
  affine, the unicity of $r$ and $g$ shows that it holds as stated; we assume
  now that $X$ and $Y$ are affine. The affine variety $X$ is defined
  over $K'$, and the domain $K'[X]$ is geometrically normal in the
  sense discussed previously.
  
  Write $Y'$ for the closure of the image of $f$. Then $Y'$ is defined
  over $K$. Moreover, if $i:Y' \to Y$ denotes the inclusion, $di_y$ is
  injective for all $y \in Y'(k)$; see e.g.
  \cite{springer-LAG}*{Exerc. 4.1.9(4)}. Since $Y'$ is again
  geometrically irreducible, we may and shall replace $Y$ by $Y'$;
  thus we assume that $f$ is a dominant morphism.  Since the tangent
  maps of $F^r_X$ are all 0, it is clear that if a suitable $r \ge 0$
  exists, it is unique.
  
  Assume that $df_x=0$ for all smooth $k$-points $x \in X(k)$; Lemma
  \ref{lem:diff-vanish} then shows that $f^*k(Y) \subset k(X)^p$.  The
  assumption on the image of $f$ means that the transcendence degree
  over $K$ of $K(Y)$ is $\ge 1$; since $k(X)$ is a finitely generated
  field extension of $k$, it follows that we may choose $r \ge 1$ such
  that $f^*k(Y) \subset k(X)^q$ for $q=p^r$ but not for $q=p^{r+1}$.

  Put $q=p^r$. We now apply Lemma \ref{lem:integral} to see that
  $f^*(K[Y]) \subset K[A^q]$.  Lemma \ref{lem:Frob-twist-q-power}
  gives then a $K$-algebra isomorphism $\phi:K[A^q] \to K[\up{X}{r}]$
  inverse to $F^r$, and we define $g:\up{X}{r} \to Y$ to have
  comorphism $\phi \circ f^*$. It is clear that $f = g \circ F^r_X$
  and that $g$ is the unique morphism with this property.
    
  The Frobenius map gives an isomorphism $F^r:k(\up{X}{r}) \to
  k(X)^q$. If $h \in K[Y]$, and if $g^*h$ is a $p$-th power in
  $k(\up{X}{r})$ then $f^*h$ is a $q'$-th power in $k(X)$, where
  $q'=p^{r+1}$.  Since $f^*k(Y)$ is not contained in $k(X)^{q'}$,
  $g^*(k(Y))$ is not contained in $k(\up{X}{r})^p$.  It then follows
  from Lemma \ref{lem:diff-vanish} that $dg_x$ is non-0 for all $x$ in
  some non-empty open subset of $X$, and the result is proved.
\end{proof}

\begin{rem}
  Let $X \subset \Aff^2$ denote the irreducible variety with
  $k$-points $\{(s,t) \mid s^p = t^p(t-1)\}$, and let $Y = \Aff^1$.
  Consider the morphism $f:X \to Y$ given on $k$-points by $f(s,t)=
  t-1$. Since $t-1=(s/t)^p$ on the open subset $U$ of $X$
  defined by $t \not = 0$, we have $df_x=0$ for each $x \in U(k)$.
  Since $X$ is over $\mathbf{F}_p$ in an obvious way, we identify $X$
  and $\up{X}{1}$; the Frobenius map $F:X \to X$ is then just $F(s,t)
  = (s^p,t^p)$.  There is a unique $\tilde g:U \to \Aff^1$ with
  $f_{\mid U} = \tilde g \circ F$; it is given on $k$-points by
  $((s,t) \mapsto s/t)$. Moreover, $d \tilde g_x \not = 0$ for each $x
  \in U(k)$.  However, there is no regular function $g$ on $X$ such
  that $g_{\mid U} = \tilde g$; thus $X$ is not normal, and the
  conclusion of the Theorem does not hold for $f$.
\end{rem}

\begin{cor}
  \label{cor:untwist}
  Let $G$ and $H$ be linear algebraic $K$-groups. Assume that $G$ is
  connected, and that $G$ is defined over the perfect subfield $K'$.
  Let $\phi:G \to H$ be a homomorphism of $K$-groups such that the
  image of $\phi$ is a positive dimensional subgroup of $H$.  There is
  a unique integer $r \ge 1$ and a unique homomorphism of $K$-groups
  $\psi:{\up{G}{r}}_{/K} \to H$ such that
  \begin{enumerate}
  \item $\phi = \psi \circ F^r_G$, and
  \item the differential $d\psi = d\psi_1$ is non-zero.
  \end{enumerate}
\end{cor}

\begin{proof}
  The $K'$-variety $G$ is geometrically irreducible; since
  $G_{/k}$ is smooth, $G$ is geometrically normal.  Hence we may apply
  Theorem \ref{theorem:untwist}; we find a unique $r \ge 0$ and a
  morphism of $K$-varieties $\psi:\up{G}{r}_{/K} \to H_{/K}$ such
  that $\psi \circ F^r_{G}$ coincides with the restriction of $\phi$
  and such that $d\psi_x$ is non-zero for $x$ in some non-empty open
  subset of $\up{G}{r}$.
  
  Since the Frobenius homomorphism $F^r_G:G \to \up{G}{r}$ is
  bijective on $k$-points, it is clear that $\psi$ is a homomorphism
  of algebraic groups. Since $d\psi_x \not = 0$ for some $x \in
  \up{G}{r}(k)$, the map induced by $\psi$ on left-invariant
  differentials in $\Omega_{\up{G}{r}/k}$ is non-0; this implies that
  $d\psi_1\not =0$ and the proof is complete.
\end{proof}

\section{Nilpotent and unipotent elements}
\label{sec:nilpotent}

We return to consideration of a strongly standard reductive $K$-group
$G$.  Let $X \in \glie$ be nilpotent. A cocharacter $\Psi:\G_m \to G$
is said to be associated with $X$ if the following conditions hold:

\begin{enumerate}
\item[(A1)] $X \in \glie(\Psi;2)$, where for any $i\in\Z$ the subspace
  $\glie(i)=\glie(\Psi;i)$ is the $i$ weight space of the torus
  $\Psi(\G_m)$ under its adjoint action on $\glie$.
\item[(A2)] There is a maximal torus $S \subset C_G(X)$ such that
  $\Psi(\G_m) \subset (L,L)$ where $L=C_G(S)$.
\end{enumerate}

With the preceding notation, $X$ is a \emph{distinguished} nilpotent
element in the Lie algebra of the Levi subgroup $L$ (see the
discussion just before Proposition \ref{prop:distinguished-parabolic}
for the definition).

If $\Psi$ is associated to $X$, the parabolic subgroup $P = P(\Psi)$
is known variously as the canonical parabolic, the Jacobson-Morozov
parabolic, or the instability parabolic (``instability flag'')
associated with $X$.  Among other things, the following result shows
 this parabolic subgroup to be independent of the choice of
cocharacter associated to $X$.

\begin{propdef}
  \label{prop:assoc-cochar}
  Let $X \in \glie(K)$ be nilpotent. 
  \begin{enumerate}
  \item There is a cocharacter $\Psi$ associated with $X$ which is
    defined over $K$.
  \item If $\Psi$ is associated to $X$ and $P=P(\Psi)$ is the
    parabolic determined by $\Psi$, then $C_G(X) \subset P$.  In
    particular, $\clie_\glie(X) \subset \Lie(P)$.
  \item Let $U$ be the unipotent radical of $C=C_G^o(X)$. Then $U$ is
    defined over $K$, and is a $K$-\emph{split} unipotent group.  If
    the cocharacter $\Psi$ is associated with $X$, then $L=C \cap
    C_G(\Psi(\G_m))$ is a Levi factor of $C$; i.e.  $L$ is connected
    and reductive, and $C$ is the semidirect product $U \cdot L$.
  \item Any two cocharacters $\Psi$ and $\Phi$ which are associated
    with $X$ are conjugate by a unique element $x \in U$. If $\Psi$ and
    $\Phi$ are each defined over $K$, then  $x \in U(K)$. 
  \item The parabolic subgroups $P(\Psi)$ for cocharacters $\Psi$
    associated with $X$ all coincide; the subgroup $P(X) = P(\Psi)$ is
    called the instability parabolic of $X$.
  \end{enumerate}
\end{propdef}

See e.g. \cite{springer-LAG}*{Chapter 14} for the notion of a
$K$-split unipotent group. We will not need to explicitly refer to
this notion here.

\begin{proof}
  The assertion (1) in the ``geometric case'' (when $K=k$) is a
  consequence of Pommerening's proof of the Bala-Carter theorem in
  good characteristic; a proof of that theorem which avoids
  case-checking has been given recently by Premet \cite{premet} using
  results in geometric invariant 
  \cite{kempf-instab}. One can deduce the assertion from Premet's
  work -- see \cite{mcninch-rat}*{Proposition 18}. Working over the
  ground field $K$, (1) was proved in \cite{mcninch-rat}*{Theorem 26}.
  
  (2) is \cite{jantzen-nil}*{Proposition 5.9}.  
  
  The first assertion of (3) is \cite{mcninch-rat}*{Theorem 28};
  notice that assumption (4.1) of \emph{loc. cit.} holds for strongly
  standard $G$, by Proposition \ref{prop:separable-orbits}. The
  semidirect product decomposition of $C$ may be found in
  \cite{jantzen-nil}*{Prop.  5.10 and 5.11}; see also
  \cite{mcninch-rat}*{Corollary 29}.
  
  We now prove (4). By (3), $C=C_G^o(X)$ is the semidirect product
  $C=U \cdot L$ of its unipotent radical $U$ and the Levi factor $L=C
  \cap C_G(\Psi(\G_m))$.  One knows by \cite{jantzen-nil}*{Lemma 5.3}
  that $\Phi = \Int(g) \circ \Psi$ for an element $g \in C$. Write $g =
  x \cdot y$ with $x \in U$ and $y \in L$.  Since $y$ centralizes
  $\Psi$, one sees that $\Phi = \Int(x) \circ \Psi$ as well. Since $U
  \cap L = \{1\}$, we see that $\Phi$ and $\Psi$ are indeed conjugate
  by the \emph{unique} element $x \in U$.
  
  Assume that $\Psi$ and $\Phi$ are defined over $K$, and write
  $S=\Psi(\G_m)$ and $S'=\Phi(\G_m)$; thus $S,S' \le C$ are tori
  defined over $K$.  We have just seen that the transporter
  \begin{equation*}
    N_C(S,S') = \{g \in C \mid gSg^{-1} = S'\}  
  \end{equation*}
  is non-empty (it has geometric points); it follows from
  \cite{springer-LAG}*{13.3.1} that $N_C(S,S')$ is defined over $K$.
  
  Choose a separable closure $\Ksep \subset k$ of the ground field
  $K$; \cite{springer-LAG}*{Theorem 11.2.7} shows that
  $N_C(S,S')(K_\sep)$ is dense in $N_C(S,S')$; we may thus find $g \in
  N_C(S,S')(K_\sep)$.  Since $S$ and $S'$ are one dimensional, and
  since $\Int(g)$ induces an isomorphism between the respective groups
  of cocharacters of these tori, we must have $\Int(g) \circ \Psi =
  \pm \Phi$. Since $g \in C$, the cocharacter $\Int(g) \circ \Psi$ is
  associated with $X$; it follows that $\Int(g) \circ \Psi = \Phi$
  e.g.  since $X \in \glie(\Int(g) \circ \Psi,2)$.

  Writing $g = y \cdot x$ with $x \in U$ and $y \in L$, we have $y =
  \lim_{t \to 0} \Psi(t)g\Psi(t^{-1})$.  By Remark
  \ref{rem:lim-over-K}, $y \in C(\Ksep)$, so that $x = y^{-1}g \in
  U(\Ksep)$. Thus $x \in U(\Ksep)$ is the unique element of $U$ for
  which $\Int(x) \circ \Psi = \Phi$. Let $\Gamma = \Gal(\Ksep/K)$ be the
  Galois group. Since $\Psi$ and $\Phi$ are $\Gamma$-stable, if
  $\gamma \in \Gamma$, we see that
  \begin{equation*}
    \Int(\gamma(x)) \circ \Psi = \Phi;
  \end{equation*}
  the unicity of $x$ shows that $x = \gamma(x)$ and we deduce that $x
  \in U(K)$ as required.
  
  To see (5), let $\Psi$ and $\Phi$ be cocharacters associated with
  $X$. Since we have $U \le C \le P(\Psi)$ by (2), it follows from (4)
  that the parabolic subgroups $P(\Psi)$ and $P(\Phi)$ are equal.
\end{proof}

Recall that a nilpotent element $X \in \glie$ is said to be
\emph{distinguished} if the connected center of $G$ is a maximal torus
of $C_G(X)$.  A parabolic subgroup $P \le G$ is said to be
distinguished if
\begin{equation*}
  \dim P/U = \dim U/(U,U) + \dim Z
\end{equation*}
where $U$ is the unipotent radical of $P$, and $Z$ is the center of
$G$.

\begin{prop}
  \label{prop:distinguished-parabolic}
  Assume that $X \in \glie$ is a distinguished nilpotent
  element. Then the instability parabolic $P=P(X)$ is a distinguished
  parabolic subgroup, and $X$ lies in the dense (Richardson) orbit
  of $P$ on $\Lie(R_uP)$. 
\end{prop}

\begin{proof}
  \cite{mcninch-rat}*{Proposition 16}.
\end{proof}

\begin{rem}
  Fixing an equivariant isomorphism $\Lambda:\UU \to \NN$ defined over
  $K$, we may say that a cocharacter $\Psi$ is associated with the
  unipotent element $u \in G$ if it is associated with $\Lambda(u)$.
  The analogous assertions of the proposition then hold for unipotent
  elements of $G$. Note that, with this definition, the notion of
  cocharacter associated with a unipotent element \emph{depends on the
    choice of $\Lambda$.}  If $\Psi$ is a cocharacter associated with
  $X=\Lambda(u)$ and if $\Lambda'$ is a second Springer isomorphism,
  easy examples show that $\Lambda'(u)$ need not be a weight vector
  for $\Psi$.  On the other hand, if $\Psi'$ is associated with
  $X'=\Lambda'(u)$, then $P(\Psi) = P(\Psi')$.  To see this, note that
  $X$ and $X'$ have the same centralizer. Fix a maximal torus $S$ of
  this centralizer and write $L=C_G(S)$; since both $\Lambda$ and
  $\Lambda'$ restrict to isomorphisms $\UU_L \to \NN_L$ (see Remark
  \ref{rem:springer-iso-restriction}), we may as well suppose that $X$
  and $X'$ are distinguished.  Since e.g.  $\Lambda'$ restricts to an
  isomorphism $U \to \Lie(U)$ where $U=R_u(P(\Psi))$, it follows that
  $X$ and $X'$ are both Richardson elements for $P(\Psi)$. Thus $\Psi$
  and $\Psi'$ are conjugate by an element of $P(\Psi)$ and it is then
  clear that $P(\Psi)=P(\Psi')$. In fact, it is even clear that $\Psi$
  and $\Psi'$ are conjugate by an element of the unipotent radical of
  $P(\Psi)$; this shows that $\Psi$ is an \emph{optimal cocharacter}
  for $X'$ (in the sense of  \cite{kempf-instab}) even though it need
  not be associated to $X'$.
\end{rem}

\section{The order formula and a generalization}

Throughout this section, $G$ is a strongly standard reductive
$k$-group defined over $K$.  Let $P$ be a parabolic subgroup of $G$;
we may fix representatives $u \in U=R_u(P)$ and $X \in \Lie(U)$ for
the dense (Richardson) $P$-orbits on $U$ and $\Lie(U)$.

Recall that if the nilpotence class of $U$ is $<p$, then
$\Lie(U)$ may be regarded as an algebraic $K$-group using the
Hausdorff formula; cf. \cite{seitz}*{\S 5}.

\begin{prop}
  \label{prop:order-result}
  Assume that $P$ is a \emph{distinguished} parabolic subgroup.
  The following conditions are equivalent:
    \begin{enumerate}
    \item $u$ has order $p$,
    \item $X^{[p]} = 0$,
    \item $\glie(\Psi;i)=0$ for all $i \ge 2p$ and some (any) cocharacter
      $\Psi$ associated to $u$ or to $X$,
    \item the nilpotence class of $U$ is $<p$.
    \end{enumerate}
\end{prop}

\begin{proof}
  The equivalence of (1) and (2) follows e.g. from
  \cite{mcninch-sub-principal}*{Theorem 35}. The equivalence of (2),
  (3) and (4) is \cite{mcninch-abelian}*{Theorem 5.4} -- note that
  there is a mis-statement (``off by 1 glitch'') concerning the
  nilpotence class in \cite{mcninch-abelian} which is explained and
  corrected in the footnote to \cite{mcninch-sub-principal}*{Lemma 11}.
\end{proof}

\begin{rem}
  Let $X$ be a distinguished nilpotent element with $X^{[p]}=0$, and
  let $U$ be the unipotent radical of the instability parabolic of
  $X$.  The proposition shows that the nilpotence class of $U$ $<p$.
  This is not true in general for nilpotent elements which are not
  distinguished.  For example, let $G = \GL_5$, and let $X \in \glie$
  be a nilpotent element with partition $(3,2)$. Then $X$ is
  distinguished in $\Lie(L)$, where $L$ is a Levi subgroup whose
  derived group is $\SL_3 \times \SL_2$. If $\Psi \in X_*(L)$ is
  associated to $X$, then $P_G(\Psi)$ is a Borel subgroup of $G$. In
  particular, if $p=3$, $X^{[p]}=0$ but a Richardson element $Y$ for
  $P_G(\Psi)$ has $Y^{[p]} \not = 0$.
\end{rem}

\begin{prop} 
  \label{prop:unipotent-exp-iso}
  Let $P$ be a \emph{distinguished} parabolic subgroup.
  If the equivalent conditions of Proposition \ref{prop:order-result}
  hold, and if $P$ is defined over $K$, then:
    \begin{enumerate}
    \item  there is a unique
      $P$-equivariant isomorphism of algebraic groups
      \begin{equation*}
        \e:\Lie(U) \to U
      \end{equation*}
      such that $d\e_0:\Lie(U) \to \Lie(U)$ is the identity.
    \item  $\e$ is defined
      over $K$.
    \item Any homomorphism $\G_a \to U$ over $K$ has the form
      \begin{equation*}
        s \mapsto \e(sX_0) \cdot \e(s^pX_1) \cdot \e(s^{p^2}X_2) \cdots
        \e(s^{p^n}X_n)
      \end{equation*}
      for some elements $X_0,X_1,\dots,X_n \in \Lie(U)(K)$ with
      $[X_i,X_j]=0$ for all $0 \le i,j \le n$.
    \end{enumerate}
\end{prop}

\begin{proof}
  Since the conditions of Proposition \ref{prop:order-result} hold,
  the unipotent radical $U = R_uP$ has nilpotence class $<p$.  In \S 5
  of \cite{seitz} -- a section contributed by J-P. Serre -- one now
  finds the necessary results. (1) and (2) follow from Proposition 5.3
  of \emph{loc. cit.}, while (3) is Proposition 5.4 of
  \emph{loc. cit.}
\end{proof}

\begin{rem}
  Recall from Remark \ref{rem:springer-iso-restriction} that the
  restriction of \emph{any} Springer isomorphism $\NN \to \UU$ gives a
  $P$-equivariant isomorphism $\Lie(U) \to U$. If $p \ge h$, there is
  always a Springer isomorphism whose restriction is $\e$. It does not
  seem to be clear (to the author, at least) whether a suitable
  analogue of this statement is true if one weakens the assumption on $p$.
\end{rem}

Recall that we may regard $G_{/k}$ as arising by base change from a
split reductive group scheme $G_{/\Z}$ over $\Z$. Write
$T_{/\Z}$ for a split maximal torus of $G_{/\Z}$.

\begin{lem}
  \label{lem:over-Z}
  Let $X \in \glie$, let $L$ be a Levi subgroup of $G$ with $X \in
  \Lie(L)$ distinguished, and let $\Psi \in X_*(L)$ be associated with
  $X$. We may find a number field $F \supset \Q$, a valuation ring
  $\Lambda \subset F$ whose residue field embeds in $k$, a standard
  Levi subgroup $M_{/\Z}$ of $G_{/\Z}$, a cocharacter $\Psi' \in
  X_*(T_{/\Z})$, and an element $Y_\Lambda \in
  \Lie(M_{/\Lambda})(\Psi';2)$ such that $(Y,M,\Psi') = g.(X,L,\Psi)$
  for some $g \in G$, where $Y=Y_\Lambda\tensor 1_k$.  Moreover, we
  may arrange that $Y_F = Y_\Lambda \tensor 1_F$ is also a Richardson
  element for the parabolic subgroup $P_{M_{/F}}(\Psi') \le M_{/F}.$
\end{lem}

\begin{proof}
  $L$ is evidently conjugate to some standard Levi subgroup $M$, which
  we may regard as arising from the Levi subgroup scheme $M_{/\Z}$.
  Replacing $X$, $L$, and $\Psi$ by a $G$-conjugate we may thus
  supposed that $L$ is standard. Replacing $(X,L,\Psi)$ by an
  $L$-conjugate, we may then assume that $X$ is a Richardson element
  for a standard distinguished parabolic of $L$.  The remainder of the
  lemma is now essentially the content of
  \cite{mcninch-abelian}*{Lemma 5.2}.
\end{proof}

\begin{prop}[Spaltenstein]
  \label{prop:spaltenstein}
  Let $\Lambda \subset F$ be a valuation ring in a number field, as in the
  previous Lemma. Let $\Psi \in X_*(T_{/\Lambda})$, let $X_\Lambda \in
  \glie_{/\Lambda}(\Psi;2)$, and assume that $\Psi$ is associated to
  $X_k$ and to $X_F$. Then
  \begin{equation*}
    \dim \clie_\glie(X_k) = \dim \clie_{\glie_{/F}}(X_F).
  \end{equation*}
\end{prop}

\begin{proof}
  This is essentially \cite{mcninch-abelian}*{Proposition 5.2} when
  $G$ is semisimple in very good characteristic.  As observed in
  \emph{loc. cit.}, it was proved by Spaltenstein for such $G$. A look
  at the proof of Spaltenstein in \cite{spaltenstein} shows that the
  result remains valid for strongly standard reductive groups [the
  only conditions on $G$ used in the proof in \cite{spaltenstein} are:
  the validity of the Bala-Carter theorem and the separability of
  nilpotent orbits].
\end{proof}

\begin{prop}
  \label{prop:p-nilpotent-cocharacter}
  Let $X \in \glie$ satisfy $X^{[p]}=0$. If $\Psi$ is a cocharacter
  associated with $X$ and if $\glie(\Psi;n) \not = 0$, then $-2p+2 \le
  n \le 2p-2$.
\end{prop}

\begin{rem}
  The analogue of the proposition for unipotent elements of order $p$
  was essentially observed by G. Seitz \cite{seitz} and is crucial to
  the proof of the existence of good $A_1$-subgroups in \emph{loc.
    cit.} It is proved for the classical groups in \cite{seitz}*{Prop.
    4.1}, and for the exceptional groups it is observed in the proof
  of \cite{seitz}*{Prop. 4.2} that it follows either from an explicit
  calculation with the associated cocharacter (``labeled diagram'')
  of each nilpotent orbit, or from some computer calculations of R.
  Lawther.
\end{rem}

\begin{proof}
  It is enough to verify the proposition for a $G$-conjugate of $\Psi$
  and $X$.  Lemma \ref{lem:over-Z} shows that, after replacing the data
  $X,L,\Psi$ by a $G$-conjugate, we may assume, as in that lemma, that
  $X$, $L$, and $\Psi$ are ``defined over $\Lambda$'' for a suitable
  valuation ring $\Lambda$. We write $X_\Lambda$ for the element of
  $\glie_{/\Lambda}$ giving rise to $X_k=X$ by base change, and we write
  $X_F = X_\Lambda \tensor 1_F \in \glie_{/F}$; note that $\Psi$ is a
  cocharacter both of $G_{/F}$ and of $G_{/k}$, and $\Psi$ is
  associated to both $X$ and $X_F$.
  
  We now contend that if $\glie(\Psi;n) \not = 0$ for some $n \ge
  2p-1$, then $\ad(X_k)^p \not = 0$; this implies the proposition. The
  proof is essentially like that of \cite{mcninch-abelian}*{Theorem
    5.4} except that we must also deal with the fact that the (in
  general, not distinguished) orbit of $X$ may not be ``even''.

  Let $\LL = \bigoplus_{i \ge -1} \glie_{/\Lambda}(\Psi;i)$, and
  $\LL^+ = \bigoplus_{i \ge 1} \glie_{/\Lambda}(\Psi;i)$.  Since we
  may embed $X_F$ in an $\lie{sl}_2(F)$-triple normalized by the image
  of $\Psi$, the representation theory of $\lie{sl}_2(F)$ implies that
  $\ad(X_F):\LL_F \to \LL^+_F$ is surjective, where the subscript
  indicates ``base change'' -- e.g.  $\LL_F = \LL \tensor_\Lambda F$.
  In view of Proposition \ref{prop:spaltenstein} and Proposition
  \ref{prop:assoc-cochar}, one knows that the kernels of the maps
  $\ad(X_k):\LL_k \to \LL_k^+$ and $\ad(X_F):\LL_F \to \LL_F^+$ have
  the same dimension. We may therefore argue as in
  \cite{mcninch-abelian}*{Proposition 5.1} and see that
  $\ad(X_k):\LL_k \to \LL_k^+$ is also surjective, hence that
  $\ad(X_k)^{n/2} \not = 0$ if $n$ is even, and that
  $\ad(X_k)^{(n+1)/2} \not = 0$ if $n$ is odd, whence our claim and
  the proposition.
\end{proof}

\section{Optimal $\SL_2$-homomorphisms.}

Throughout this section, $G$ will denote a strongly standard reductive
$K$-group.  We first ask the reader's patience while we fix some
convenient notation for $\SL_2$. We choose the standard basis for
$\lie{sl}_2$:
\begin{equation*}
  X_1 = 
  \begin{pmatrix}
    0 & 1 \\
    0 & 0
  \end{pmatrix},
  \quad  
  H_1 = 
  \begin{pmatrix}
    1 & 0 \\
    0 & -1
  \end{pmatrix},
  \quad \text{and} \quad
  Y_1 =  \begin{pmatrix}
    0 & 0 \\
    1 & 0
  \end{pmatrix}. 
\end{equation*}  
Now put:
\begin{equation*}
  x_1(t) = 
  \begin{pmatrix}
    1 & t \\
    0 & 1
  \end{pmatrix} 
  \quad \text{and} \quad
  y_1(t) = 
  \begin{pmatrix}
    1 & 0 \\
    t & 1
  \end{pmatrix} \text{\ for\ } t \in k,
\end{equation*}
and write
$\X = \{ x_1(t) \mid t\in k\}$ and $\X^- = \{y_1(t) \mid t \in k\}$.
Finally, write
\begin{equation*}
  \T = \left \{ 
    \begin{pmatrix}
      t & 0\\
      0 & t^{-1}
    \end{pmatrix}\mid t \in k^\times
  \right\}
\end{equation*}
for the standard maximal torus of $\SL_2$. 

We fix once and for all one of the two isomorphisms $\G_m \iso \T$, so
that if $\phi:\SL_2 \to G$ is a homomorphism, it determines a
cocharacter $\Psi= \phi_{\mid \T} \in X_*(G)$ by restriction to $\T$; explicitly,
$\Psi$ is given by the rule $$\Psi(t) = \phi(
\begin{pmatrix}
  t & 0 \\
  0 & t^{-1}
\end{pmatrix}) \quad \text{for $t \in k^\times$}.$$

\begin{defin}
  The homomorphism $\phi:\SL_2 \to G$ is an \emph{optimal
    $\SL_2$-homomorphism} if the cocharacter $\Psi =
  \phi_{\mid \T}$ is associated to the nilpotent element $X =
  d\phi(X_1) \in \glie$.  Briefly, we say that $\phi$ is optimal for
  $X$.
\end{defin}

We first recall that the main result of \cite{mcninch-sub-principal}
shows that optimal homomorphisms always exist.  More precisely, let $X
\in \glie$ with $X^{[p]}=0$, and let $\Psi$ be a cocharacter
associated with $X$. If $S$ is a maximal torus of $C_\Psi$, then $X$
is distinguished in $\Lie(L)$ where $L=C_G(S)$.  We may apply
Proposition \ref{prop:unipotent-exp-iso} to $P_L(\Psi)$; let
$\e:\Lie(U) \to U$ be the isomorphism of that proposition, where we
have written $U$ for the unipotent radical of $P_L(\Psi)$.  Now the
main result of \cite{mcninch-sub-principal} says the following:
\begin{prop}
  \label{prop:optimal-exist}
  There is an optimal $\SL_2$-homomorphism $\phi$ for $X$ with the
  following properties: 
  \begin{enumerate}
  \item $\phi_{\mid \T} = \Psi$, and
  \item $\phi(x(t)) = \e(tX)$ for each $t \in k$.
  \end{enumerate}
\end{prop}

We wish to see that $\e(tX)$ is independent of the choice of the
maximal torus $S$ of $C_\Psi$. For this, we will use the following
result due to Seitz; the result is essentially \cite{seitz}*{Prop. 4.2}.

\begin{prop}[Seitz]
  \label{prop:tilting}
  Let $\Lambda \subset F$ be a valuation ring in a number field whose
  residue field is embedded in $k$, let $\LL$ be a $\Lambda$ lattice,
  and let $\rho_{/\Lambda}:\SL_{2/\Lambda} \to \GL(\LL)$ be a
  representation over $\Lambda$. Assume that
  \begin{enumerate}
  \item all weights of the standard maximal $\Lambda$-torus
    $\T_{/\Lambda}$ on $\LL$ are $\le 2p-2$,
  \item the representation $\rho_{/k}$
    of $\SL_{2/k}$ is self-dual,
%the representations $\rho_{/F}$ of $\SL_{2/F}$ and $\rho_{/k}$
%    of $\SL_{2/k}$ are self-dual,
  \item the dimension of the fixed point space of $u_F=\rho_{/F}(
    \begin{pmatrix}
      1 & 1\\
      0 & 1
    \end{pmatrix})$ on $\LL_F$ is  the same as the dimension of the fixed
    point space of $u_k=\rho_{/k} (   \begin{pmatrix}
      1 & 1\\
      0 & 1
    \end{pmatrix})$ on $\LL_k$.
  \end{enumerate}
  Then the representation $(\rho_{/k},\LL_k)$ is a \emph{tilting
    module} for $\SL_{2/k}$.
\end{prop}

\begin{proof}
  One decomposes the $\SL_{2/k}$-module $\LL_k$ according to the
  blocks of $\SL_{2/k}$. In view of the assumption on the weights of
  $\T_{/k}$ on $\LL_k$, the blocks that can conceivably occur are
  those of the simple modules $L(d)$ with $0 \le d < p$. The summand
  corresponding to the block for $d=p-1$ is isomorphic to
  $L(d)^{v(d)}$ for some integer $v(d) \ge 0$. Otherwise, the summand
  corresponding to a block with $d < p-1$ is isomorphic to a module of
  the form
  \begin{equation*}
    T(c_d)^{r(d)} \oplus W(c_d)^{s(d)} \oplus  (W(c_d)^\vee)^{t(d)} 
    \oplus L(c_d)^{u(d)}  \oplus L(d)^{v(d)}
  \end{equation*}
  where $c_d = 2p-2-d$ and where the exponents
  $r(d),s(d),t(d),u(d),v(d)$ are non-negative integers.  [We are using
  Seitz's notation for
  $\SL_{2/k}$-representations: $W(d)$ is the Weyl module with high
  weight $d$, and $T(d)$ is the indecomposable tilting module with
  high weight $d$; cf. \cite{seitz}*{\S 2}.]
  
  The assumption (2) implies that $s(d) = t(d)$ for all $0 \le d <
  p-1$.  As in \cite{seitz}*{Prop. 4.2}, one now expresses the
  dimensions of the fixed point spaces of $u_k$ and $u_F$ in terms of
  the exponents and finds that $u(d) = s(d) = t(d) = 0$ for all $d$.
  Thus $\LL_k$ is the direct sum of various simple tilting modules $L(d)$ for
  $0 \le d<p$, and various indecomposable tilting modules $T(c_d) =
  T(2p-2-d)$ for $0 \le d < p-1$, so indeed $\LL_k$ is a tilting module.
\end{proof}

\begin{prop}
  \label{prop:exp-centralizer}
  With notation as above, we have
  \begin{enumerate}
  \item $C_G^o(X) = C_G^o(\e(X))$; in particular, $\Psi(\G_m)$
    normalizes $C_G^o(\e(X))$.
  \item $C_G^o(\e(X)) = C_G^o(\e(tX))$ for each $t \in k^\times$.
  \end{enumerate}
\end{prop}

\begin{proof}
  If $X$ is distinguished, (1) holds since $\e$ is $P=P(\Psi)$
  equivariant, since $\e(X) \in R_u(P)$ is again a Richardson element,
  and since $C_G(X),C_G(\e(X)) \le P$ by Proposition
  \ref{prop:assoc-cochar}. [In fact, $C_G(X) = C_G(\e(X))$ always
  holds in this case.] It remains to prove (1) when $X$ is no longer
  distinguished; we essentially follow the proof in
  \cite{seitz}*{Lemma 6.3}.
  
  By the unicity of $\e$, it is enough to prove the result with $L$,
  $\Psi$, and $X$ replaced by a $G$-conjugate.  We will regard
  $G=G_{/k}$ as arising by base change from the split reductive group
  scheme $G_{/\Z}$ over $\Z$; let $T_{/\Z}$ be a $\Z$-split maximal
  torus of $G_{/\Z}$.
  
  According to Lemma \ref{lem:over-Z}, we may find a suitable
  valuation ring in a number field $\Lambda \subset F$ and assume that
  the Levi subgroup $L$ contains $T_{/k}$ and arises by base change
  from a standard split reductive Levi subgroup scheme $L_{/\Z} \le
  G_{/\Z}$ containing $T_{/\Z}$, that $\Psi \in X_*(T_{/\Z})$, and
  that the nilpotent element $X_\Lambda \in
  \Lie(L_{/\Lambda})(\Psi;2)$ gives $X$ on base change.

  After possibly enlarging $\Lambda$ and $F$,
  \cite{mcninch-sub-principal}*{Theorem 13} gives a homomorphism 
  \begin{equation*}
    f:\SL_{2/\Lambda}
  \to G_{/\Lambda}
  \end{equation*}
  such that the restriction of $f$ to the subgroup
  scheme $
  \begin{pmatrix}
    1 & * \\
    0 & 1
  \end{pmatrix}$ of $\SL_{2/\Lambda}$
  is given by $t \mapsto \e(tX_\Lambda)$, where $X_\Lambda \in
  \glie_{/\Lambda}$ gives $X$ upon extension of scalars to $k$ (recall
  from \cite{seitz}*{Prop.  5.1} that $\e$ is indeed defined over
  $\Z_{(p)}$ hence over $\Lambda$). Moreover, the restriction of $f$
  to the standard maximal torus of $\SL_{2/\Lambda}$ gives the
  cocharacter $\Psi$ of $T_{/\Lambda}$.
  
  Since $G$ is strongly standard, its adjoint representation is
  self-dual.  Together with Proposition \ref{prop:spaltenstein}, this
  shows that we may apply Proposition \ref{prop:tilting} to the
  representation $\Ad \circ f:\SL_{2/\Lambda} \to
  \GL(\Lie(G_{/\Lambda}))$.  Thus the $\SL_2$-representation $(\Ad
  \circ f_{/k},\glie)$ is a tilting module, and it follows from
  \cite{seitz}*{Lemma 2.3(d)} that
  \begin{equation*}
    \clie_\glie(\e(tX)) = \clie_\glie(X)
  \end{equation*}
  for each $t \in k^\times$.  The orbits of $\e(tX)$ and $X$ are
  separable by Proposition \ref{prop:separable-orbits}; thus we know
  that $\Lie C_G(\e(tX)) = \Lie C_G(X)$. In particular, $C_G(X)$ and
  $C_G(\e(X))$ have the same dimension; assertion (1) will follow if we show
  that $C_G^o(X) \le C_G^o(\e(X))$.

  For any connected linear group $H$, we write $H_t$ for the subgroup
  generated by the maximal tori in $H$. Applying
  \cite{springer-LAG}*{13.3.12}, to the group $H=C_G^o(X)$, we find
  that $H$ is generated by $H_t$ and $C_H(S)$, where $S$ is our fixed maximal
  torus of $H$; i.e.
  \begin{equation}
    \label{eq:generate-centralizer}
    H=\langle H_t,C_H(S) \rangle.
  \end{equation}
  
  Working for the moment inside the Levi subgroup $L=C_G(S)$ of $G$,
  the ``distinguished'' case of part (1) of the proposition means that
  $C_H(S)=C_L(X)=C_L(\e(X))$; in particular $C_H(S)$ centralizes
  $\e(X)$. So according to \eqref{eq:generate-centralizer}, the
  containment $H \le C_G^o(\e(X))$, and hence (1), will follow if we
  just show that $\e(X)$ is centralized by each maximal torus $T$ of
  $C_G(X)$.  Since $\clie_\glie(\e(X)) = \clie_\glie(X) = \Lie
  C_G(X)$, one knows that $\e(X)$ centralizes $\Lie(T)$. We claim that
  $(*) \ C_G(T) = C_G(\Lie(T))$; this shows that $T$
  centralizes $\e(X)$ as desired.
  
  Write $M=C_G(T)$.  Since $T$ is a maximal torus of $C_G^o(X)$, it
  follows that $T$ is a maximal torus of the center of $M$.  Thus
  $(*)$ is a consequence of the next lemma (Lemma
  \ref{lem:torus-centralizer}), and (1) is proved.  For (2), notice
  that if $s^2 =t$, we have by (1) that
  \begin{equation*}
    C_G^o(\e(X)) = \Psi(s)C_G^o(\e(X))\Psi(s^{-1}) = 
    C_G^o(\e(\Ad(\Psi(s))X)) = C_G^o(\e(tX)).
  \end{equation*}
\end{proof}

\begin{lem}
  \label{lem:torus-centralizer}
  Let $G$ be a strongly standard reductive group, let $T \le G$ be a
  torus, and write $M=C_G(T)$. If $T$ is a maximal torus of the center
  of $M$, then $C_G(T) = C_G(\Lie(T))$.
\end{lem}

\begin{proof}
  We essentially just reproduce the proof of \cite{seitz}*{Lemma 6.2}.
  Let $T_0$ be a maximal torus of $G$ containing $T$. Denote by $R
  \subset X^*(T_0)$ the roots of $G$ and by $R_L \subset R$ the roots
  of $L$.  Choose a system $\alpha_1,\dots,\alpha_r \in X_*(T_0)$ of
  simple roots for $G$ such that $\alpha_1,\dots,\alpha_t$ is a
  system of simple roots for $M=C_G(T)$ (so $t \le r$). If we write
  $U_\alpha \le G$ for the root subgroup corresponding to $\alpha \in
  R$, then $U_\alpha \le L$ for $\alpha \in R_L$; moreover,
  \begin{equation*}
    C_G(T) = \langle T_0; U_\alpha \mid \alpha_{\mid T} = 1\rangle, 
    \ \text{and}\ 
    C_G(\Lie(T)) = \langle T_0; U_\alpha \mid 
    d\alpha_{\mid \Lie(T)} = 0\rangle.
  \end{equation*}
  We have always $C_G(T) \le C_G(\Lie(T))$. If the Lemma were not
  true, there would be some root $\beta$ of $G$ such that $\beta_{\mid
    T} \not = 1$ but $d\beta_{\mid \Lie(T)} = 0$. We may write $\beta
  = \alpha + \sum_{i=t+1}^r c_i \alpha_i$ with $\alpha \in R_L$.
  Since $p$ is good, the $c_i$ are integers with $0 \le c_i < p$
  \cite{springer-steinberg}*{I.4.3}. Since $\beta_{\mid T} \not =1$,
  it follows that $c_j$ is non-zero in $k$ for some $t+1 \le j \le r$.
  
  Since $G$ and $M$ are strongly standard,
  \cite{springer-steinberg}*{Corollary I.5.2} implies that
  $\zlie(\glie) = \Lie Z(G)$ and $\zlie(\mlie) = \Lie Z(M)$ (where
  $\zlie(?)$ denotes the center of a Lie algebra, and $Z(?)$ that of a
  group). We thus have $\dim T = \dim \zlie(\glie) + (r-t)$.  It
  follows that $\{d\alpha_{t+1},\cdots,d\alpha_r\}$ is a linearly
  independent subset of $\Lie(T)^\vee$ (the dual space of $\Lie(T)$).
  In particular, there is $A \in \Lie(T)$ such that 
  \begin{equation*}
    d\alpha_i(A) = \delta_{i,j}.
  \end{equation*}
  But then $d\beta(A) = c_j \not = 0$, contradicting the choice of
  $\beta$. This completes the proof.
\end{proof}

\begin{rem}\label{rem:exp-commute}
  If $S,S' \le C_\Psi$ are maximal tori, let us write $U$ and $U'$ for
  the unipotent radicals of the distinguished parabolic subgroups
  $P_L(\Psi) \le L$ and $P_{L'}(\Psi) \le L'$ where $L=C_G(S)$ and
  $L'=C_G(S')$. If $\e:\Lie(U) \to U$ and $\e':\Lie(U') \to U'$ are
  the isomorphisms of Proposition \ref{prop:unipotent-exp-iso}, then
  $\e(tX)=\e'(tX)$ for each $t \in k$.  Indeed, we may choose $g \in
  C_\Psi^o(X)$ with $gSg^{-1} = S'$. It is then clear that $U' =
  gUg^{-1}$ and the uniqueness statement of Proposition
  \ref{prop:unipotent-exp-iso} shows that $\e' = \Int(g) \circ \e
  \circ \Ad(g^{-1}):\Lie(U') \to U'$.  Let $t \in k^\times$.
  Proposition \ref{prop:exp-centralizer} shows that $g$ centralizes
  $\e(tX)$ in addition to $X$.  So indeed
  \begin{equation*}
    \e'(tX) = \Int(g) \circ \e \circ \Ad(g^{-1})(tX) =
    \Int(g) \circ \e(tX) = \e(tX)
  \end{equation*}
  as asserted.
\end{rem}

Now let $\phi:\G_a \to G$ be an injective homomorphism of algebraic
groups with $X = d\phi(1)$, and assume that the cocharacter $\Psi$
associated to $X$ has the property that
\begin{equation*}
  \Psi(t) \phi(s )\Psi(t^{-1}) = \phi(t^2 s) \quad \text{for each}\ 
  t \in k^\times \text{ and } s \in k.
\end{equation*}
Since $\phi$ is injective, the cocharacter $\Psi$ is non-trivial; this
means in particular that $X \ne 0$ and so  $d\phi$ is non-zero.
%
%If $X=0$ then $\Psi$ is trivial and hence $\phi = 1$.  
%

We remark that the homomorphism $h:\G_a \to G$ given by $t \mapsto
\e(tX)$ is injective.  Indeed, as in the proof of Proposition
\ref{prop:exp-centralizer}, there is an optimal homomorphism
$f:\SL_2 \to G$ such that $h(s) = f(x_1(s))$ for $s \in \G_a$. The
group $\SL_2$ is almost simple; its unique normal subgroup is
contained in each maximal torus. In particular, $\ker h$ is trivial as
asserted.

Fix now  a maximal torus $S$ of $C_G(X)$ centralized by the
image of $\Psi$, and hence a Levi subgroup $L=C_G(S)$ such that
$\Psi(\G_m) \le L$ and $X \in \Lie(L)$.

\begin{prop}
  \label{prop:additive-homom}
  With $\phi$ and $\Psi$ as above, we have $\phi(t) = \e(tX)$ for each
  $t \in k$, where $\e:\Lie(U) \to U$ is the isomorphism of
  Proposition \ref{prop:unipotent-exp-iso} for the unipotent radical
  $U$ of the distinguished parabolic subgroup $P_L(\Psi) \le L$. In
  particular, $\phi(\G_a) \le L$.
\end{prop}

\begin{proof}
  Notice that $\phi(s) \in C_G^o(X)$ for all $s \in \G_a$.  According
  to Proposition \ref{prop:exp-centralizer} this shows that $\phi(s) \in
  C_G^o(\e(tX))$ for all $t \in k^\times$, hence that
  \begin{equation*}
    s \mapsto \e(-sX) \cdot \phi(s)
  \end{equation*}
  is a homomorphism $\phi_1:\G_a \to G$. Moreover,
  $\Psi(t)\phi_1(s)\Psi(t^{-1}) = \phi_1(t^2s)$ for $t\in k^\times$
  and $s \in k$, and a quick calculation shows $d\phi_1$ to be trivial.

  Assume that the proposition is not true, hence that $\phi_1 \not =
  1$; it has positive dimensional image and so by Corollary
  \ref{cor:untwist}  there is a homomorphism $\phi_2:\G_a
  \to G$ and an integer $r \ge 1$ such that $\phi_1 = \phi_2 \circ
  F^r$, where $F$ denotes the Frobenius morphism for $\SL_2$, and such
  that $d\phi_2 \not = 0$.  On the additive group, $F$ is given by $s
  \mapsto s^p$, so we know that $\phi_1(s)=\phi_2(s^{p^r})$ for $s \in
  k$. [Notice we have used the fact that $\G_a$ is defined over
  $\F_p$, so that $\G_a$  identifies with $\up{\G_a}{r}$ for $r \ge
  0$.]
  
  Observe that if $\phi_1(s_0) = 1$ for some $s_0 \not = 0$, then $1 =
  \phi_1(s_0) = \e(-s_0X)\phi(s_0)$ so that $\e(s_0X) = \phi(s_0)$;
  applying $\Int(\Psi(t))$ for $t \in k^\times$, we see that $\e(sX) =
  \phi(s)$ for all $s \in k$, so that $\phi_1 = 1$.  Thus if $\phi_1
  \not = 1$, then $\phi_1$ is an injective map on the points of
  $\G_a$. It is then clear that $\phi_2$ is injective as well [since
  $d\phi_2$ is non-zero, $\phi_2$ is even an injective homomorphism of
  algebraic groups].

  Since $\Psi(\G_m)$ normalizes the image of $\phi_2$, we have
  $\Psi(t)\phi_2(s)\Psi(t^{-1}) = \phi_2(t^ns)$ for some $n \in \Z$.
  Let now $t \in k^\times$ and $s \in k$. Then
  \begin{equation*}
    \phi_1(t^2s) = \Psi(t)\phi_1(s)\Psi(t^{-1}) = 
    \Psi(t)\phi_2(s^{p^r})\Psi(t^{-1})
    = \phi_2(t^n s^{p^r});
  \end{equation*}
  since $\phi_1$ and $\phi_2$ are injective, we have $(t^2s)^{p^r} =
  t^ns^{p^r}$ for all $t \in k^\times$ and $s\in k$. It follows that
  $n = 2p^r$.
  
  Denoting by $0\not =Y$ an element in the image of $d\phi_2$, it is
  clear that $\Ad(\Psi(t))Y = t^{2p^r}Y$ so that $Y \in
  \glie(\Psi;2p^r)$. Since $r \ge 1$, since $\Psi$ is associated with
  $X$, and since $X^{[p]}=0$, this contradicts Proposition
  \ref{prop:p-nilpotent-cocharacter}; hence $\phi_1 =1$ and
  $\phi(s)=\e(sX)$ for all $s \in k$ as asserted.
\end{proof}

\begin{rem}
  Assume that $p \ge h$, where $h$ is the Coxeter number of $G$.  Then
  the nilpotence class of the unipotent radical $U$ of a Borel
  subgroup $B$ of $G$ is $<p$. Thus there is a $B$-equivariant
  isomorphism $\e:\Lie(U) \to U$ as in Proposition
  \ref{prop:unipotent-exp-iso}.  Fix a regular nilpotent element $X
  \in \Lie(U)$ and write $u = \e(X)$.  According to Proposition
  \ref{prop:serre-springer-iso}, there is a unique Springer
  isomorphism $\Lambda:\UU \to \NN$ with $\Lambda(u)= X$.  It is then
  clear by the unicity of $\e$ that $\Lambda^{-1}_{\mid \Lie(U)} = \e$
  for the unipotent radical $U$ of \emph{any} Borel subgroup of $G$.
  Since the unipotent radical $V$ of any parabolic subgroup $P$ of $G$
  is contained in that of some Borel subgroup, it is then clear that
  ${\Lambda^{-1}}_{\mid \Lie(V)}$ is the isomorphism of Proposition
  \ref{prop:unipotent-exp-iso} (of course, the nilpotence class of $V$
  is $<p$). This permits for these $p$ a simple proof of Proposition
  \ref{prop:exp-centralizer} and hence of Proposition
  \ref{prop:additive-homom} (i.e. a proof independent of the tilting
  module considerations of Proposition \ref{prop:tilting})
\end{rem}

\subsection{Conjugacy of optimal $\SL_2$ homomorphisms}

The goal of this paragraph is to show that any two optimal
$\SL_2$-homomorphisms for $X$ are conjugate by an element of
$C_G^o(X)$. 

Let $\phi$ be an optimal $\SL_2$-homomorphism for $X \in \glie$ with
cocharacter $\Psi=\phi_{\mid \T}$.  Choose a maximal torus $S \le
C_\Psi$, so that $X$ is distinguished in $\Lie(L)$, where $L=C_G(S)$
is a Levi subgroup of $G$. If $\phi$ is defined over $K$, then the
maximal torus $S$ -- and so also $L$ -- may be chosen over $K$.

We will write $P_L = P_L(\Psi)$ for the parabolic subgroup of $L$
determined by the cocharacter $\Psi$, and $U$ for the unipotent
radical of $P_L$. Denote by $\e:\Lie(U) \to U$ the
unique $P_L$-equivariant isomorphism of
Proposition \ref{prop:unipotent-exp-iso}. 

\begin{prop}
  \label{prop:sl2-exponential-map}
   \begin{enumerate}
   \item The torus $S$ centralizes $\phi(\X)$; in particular, $\phi(\X)
     \subset U$.
  \item $\phi(x_1(t)) = \e(tX)$ for each $t \in k$.
  \item For each $t \in k^\times$, $C_G^o(X) = C_G^o(u_t))$ where $u_t =
    \phi(x_1(t))$.
  \end{enumerate}
\end{prop}

\begin{proof}
  We apply the result of Proposition \ref{prop:additive-homom}; that
  proposition shows that $\phi(t) = \e(tX)$. (1) and (2) are then
  immediate, and (3) follows from Proposition
  \ref{prop:exp-centralizer}.
\end{proof}

\begin{prop}
  \label{prop:optimal-sl2-Levi}
  The image of $\phi$ lies in the derived group of the Levi subgroup
  $L=C_G(S)$.  
\end{prop}

\begin{proof}
  Since $\SL_2$ is equal to its own derived group, we only must see
  that the image of $\phi$ lies in $L$.

  Now write 
  \begin{equation*}
    Y  =  d\phi (Y_1) \in \glie
    \quad \text{and} \quad
    u^-_t = \phi (y_1(t)) \in G \quad \text{for $t \in k$}.
  \end{equation*}
  Since $\SL_2$ is generated by the subgroups $\X$ and $\X^-$, it
  suffices to show that $u_t, u^-_t \in L=C_G(S)$ for all $t \in
  k^\times$. Fix $t \in k^\times$.  It was proved in Proposition
  \ref{prop:sl2-exponential-map}(1) that $u_t \in L$.
  
  Now, there is $g \in \phi(\SL_2)$ with $gu_tg^{-1} = u^-_t$ and $\Ad(g)X =
  Y$.  Together with Proposition \ref{prop:sl2-exponential-map}, this
  implies that $C_G^o(u^-_t) = C_G^o(Y)$ for $t \in k^\times$. So the
  proof is complete once we show that $S \le C_G(Y)$.
  
  Since $S$ and the image of $\Psi$ commute, $\glie(\Psi;-2)$ is
  $S$-stable and is thus a direct sum of $S$-weight spaces
  \begin{equation*}
    \glie(\Psi;-2) = \sum_{\gamma \in X^*(S)} \glie(\Psi;-2)_\gamma.
  \end{equation*}
  Hence, we may write $Y \in \glie(\Psi;-2)$ as a sum of $S$-weight
  vectors:
  \begin{equation*}
    Y = \sum_\gamma Y_\gamma \quad \text{with~} 
    Y_\gamma \in \glie(\Psi;-2)_\gamma.
  \end{equation*}
  We need to show that $Y=Y_0$, or equivalently that $Y_\gamma = 0$
  for $\gamma \not = 0$.
  
  As $\Psi$ is associated to $X$, it follows from Proposition
  \ref{prop:assoc-cochar} that $\clie_\glie(X) \subseteq
  \sum_{i \ge 0} \glie(\Psi;i)$. Since $S$ centralizes $X$, it follows
  that $\ad(X):\glie(\Psi;2) \to \glie(\Psi;0)$ is an injective map
  of $S$-representations. Writing $H = d\Psi(1) \in \glie$, we have
  $\ad(X)Y = [X,Y] = H \in \glie(\Psi;0)_0$. Since $\ad(X)Y_\gamma \in
  \glie(\Psi;0)_\gamma$, the injectivity of $\ad(X)$ implies that
  $Y_\gamma = 0$ unless $\gamma = 0$, as desired. Thus $Y=Y_0$ and the
  proof is complete.
\end{proof}

\begin{prop}
  \label{prop:optimal-equality}
  Let $X \in \glie$ satisfy $X^{[p]}=0$. If $\phi_1$ and $\phi_2$ are
  optimal $\SL_2$-homomorphisms for $X$ and if ${\phi_1}_{\mid \T} =
  {\phi_2}_{\mid \T}$, then $\phi_1=\phi_2$.
\end{prop}

\begin{proof}
  Combined with Proposition \ref{prop:optimal-sl2-Levi}, the
  hypotheses yield a maximal torus $S \le C_G(X)$ such that the image
  of $\phi_i$ lies in $L = C_G(S)$ for $i=1,2$.  Thus we may replace
  $G$ by the strongly standard reductive group $L$ and so suppose that
  $X$ is \emph{distinguished}.
  
  Proposition \ref{prop:sl2-exponential-map} shows that
  $\phi_1(x_1(t)) = \e(tX) = \phi_2(x_1(t))$ for all $t \in k$.  It
  follows that $\phi_1$ and $\phi_2$ coincide on the Borel subgroup
  $B=\T \X$ of $\SL_2$. Using this, we argue that $\phi_1$ and
  $\phi_2$ coincide on all of $\SL_2$. Indeed, consider the morphism
  of varieties $\SL_2 \to G$ given by
  \begin{equation*}
    g \mapsto \phi_1(g)\phi_2(g^{-1}).
  \end{equation*}
  Since the $\phi_i$ are homomorphisms, this morphism factors through
  the flag variety $\SL_2/B = \mathbf{P}^1$ (the projective line);
  since $\mathbf{P}^1$ is an irreducible complete variety, and since
  $G$ is affine, this morphism must be constant. The proof is
  complete. 
  %\footnote{I originally used a more complicated ``big
%    cell'' argument to prove this proposition; cf.
%    \cite{mcninch-sub-principal}*{5.2}. Both Jens C. Jantzen and the
%    referee suggested the simpler argument used here.}
\end{proof}

\begin{cor}
  \label{cor:optimal-centralizer}
  If $\phi$ is an optimal homomorphism, let as usual $X = d\phi(X_1)$
  and $\Psi = \phi_{\mid \T}$. Then the centralizer of $\phi(\SL_2)$
  is $C_\Psi = C_G(X) \cap C_G(\Psi(\G_m))$.
\end{cor}

\begin{proof}
  This is just a restatement of the previous proposition.
\end{proof}

\begin{theorem}
  \label{theorem:optimal-conjugate}
  Suppose that $G$ is strongly standard, and that $X \in \glie$
  satisfies $X^{[p]}=0$.  Then any two optimal $\SL_2$-homomorphisms
  for $X$ are conjugate by a unique element of the unipotent radical of
  $C_G^o(X)$.
\end{theorem}
  
\begin{proof}
  Let $\phi_1,\phi_2$ be optimal $\SL_2$-homomorphisms for $X$, and
  write $\Psi_i = {\phi_i}_{\mid \T}$ for the corresponding
  cocharacters. According to Proposition \ref{prop:assoc-cochar}, the
  cocharacters $\Psi_1$ and $\Psi_2$ associated with $X$ are conjugate
  by a unique element of the unipotent radical $U$ of $C_G^o(X)$.  Replacing
  $\phi_2$ by a $U$-conjugate, we may thus suppose that
  $\Psi_1=\Psi_2$.  It then follows from Proposition
  \ref{prop:optimal-equality} that $\phi_1=\phi_2$.
\end{proof}

\subsection{Uniqueness of a principal homomorphism}

Suppose that $X \in \glie$ is a \emph{distinguished} nilpotent
element.  Then any cocharacter $\Psi \in X_*(G)$ with $X \in
\glie(\Psi;2)$ is associated to $X$.  In particular, if $\phi:\SL_2
\to G$ is any homomorphism with $d\phi(X_1) = X$, then $\Psi =
\phi_{\mid \T}$ is a cocharacter associated with $X$; thus $\phi$ is
optimal.

An application of Theorem \ref{theorem:optimal-conjugate} now gives:
\begin{prop}
  If $\phi_1,\phi_2:\SL_2 \to G$ are homomorphisms such that
  $d\phi_1(X_1) = d\phi_2(X_1)=X$ is a distinguished nilpotent
  element, then $\phi_1$ and $\phi_2$ are conjugate by an element of
  $C_G^o(X)$.
\end{prop}

A \emph{principal homomorphism} $\phi:\SL_2 \to G$ is one
for which $d\phi(X_1)$  is a regular nilpotent element.
Since a regular nilpotent element is distinguished, we have:
\begin{prop}
  A principal homomorphism is optimal. Any two principal homomorphisms
  are conjugate in $G$.
\end{prop}

\subsection{Optimal homomorphisms over ground fields}
Recall that $K$ is an arbitrary ground field. The following theorem gives
both an existence result and a conjugacy result for optimal
homomorphisms over the ground field $K$. If $X \in \glie(K)$,
write $C = C_G^o(X)$ for its connected centralizer; recall by
Proposition \ref{prop:assoc-cochar} that the unipotent radical
of $C$ is defined over $K$.
\begin{theorem}
  \label{theorem:optimal-over-ground-field}
  Let $G$ be a strongly standard reductive $K$-group, and let $X \in
  \glie(K)$ satisfy $X^{[p]}=0$. 
  \begin{enumerate}
  \item There is an optimal $\SL_2$-homomorphism $\phi$ for $X$ which
    is defined over $K$.
  \item Let $U$ be the unipotent radical of $C = C_G^o(X)$. Any two
    optimal $\SL_2$-homomorphism for $X$ defined over $K$ are
    conjugate by a unique element of $U(K)$.
  \end{enumerate}
\end{theorem}

\begin{proof}
  To prove (1), we need first to quote a more precise form of
  Proposition \ref{prop:optimal-exist}. The proof of that Proposition
  given in \cite{mcninch-sub-principal} shows that there is a
  nilpotent element $X''$ in the orbit of $X$ which is rational over
  the separable closure $K_\sep$ of $K$ in $k$ and an optimal
  $\SL_2$-homomorphism $\phi''$ for $X''$ defined over $K_\sep$.
  Since the orbit of $X$ is separable, one can mimic the proof of
  \cite{springer-LAG}*{12.1.4} to see that $X$ and $X''$ are conjugate
  by an element rational over $K_\sep$. Indeed, let $\mathcal{O}$ be
  the orbit of $X$ and let $\mu:G \to \mathcal{O}$ be the orbit map
  $\mu(g) = \Ad(g)X$.  The separability of the orbit $\mathcal{O}$
  means that $d\mu_1:T_1(G) \to T_X(\mathcal{O})$ is surjective, and
  it follows for each $g \in G$ that $d\mu_g:T_g(G) \to T_{\Ad(g)X}
  (\mathcal{O})$ is surjective.  It follows from
  \cite{springer-LAG}*{11.2.14} that the fiber $\mu^{-1}(X'')$ is
  defined over $K_\sep$, so that by \cite{springer-LAG}*{11.2.7} there
  is a $K_\sep$-rational point $g$ in this fiber. It follows that
  $\phi'=\Int(g) \circ \phi''$ is an optimal $\SL_2$-homomorphism for
  $X$ which is defined over $K_\sep$.
  
  According to Proposition \ref{prop:assoc-cochar}, we can find a
  cocharacter $\Psi$ associated with $X$ which is defined over $K$.
  Writing $C=C_G^o(X)$, that same Proposition shows that the
  cocharacters $\Psi$ and $\Psi'=\phi'_{\mid \T}$ are conjugate by an
  element $h \in C(\Ksep)$ [in fact, $h$ can be chosen to be a
  $\Ksep$-rational element of the unipotent radical of $C$].
  
  It now follows that $\phi=\Int(h^{-1}) \circ \phi'$ is an optimal
  $\SL_2$-homomorphism for $X$ which is defined over $K_\sep$. We
  argue that $\phi$ is actually defined over $K$. Let $\gamma \in
  \Gal(K_\sep,K)$.  Then $\phi_\gamma = \gamma \circ \phi \circ
  \gamma^{-1}:\SL_2 \to G$ is another optimal $\SL_2$-homomorphism for
  $X$; since $\Psi = \phi_{\mid \T}$ is defined over $K$, $\phi_{\mid
    \T} = {\phi_\gamma}_{\mid \T}$.  Thus Proposition
  \ref{prop:optimal-equality} shows that $\phi=\phi_\gamma$. Since
  $\phi$ is defined over $K_\sep$, Galois descent (e.g.
  \cite{springer-LAG}*{Cor. 11.2.9})  shows that $\phi$ is defined
  over $K$.
  
  We now give the proof of (2), which is the same as the proof of
  Theorem \ref{theorem:optimal-conjugate}. If $\phi$ and $\psi$ are
  optimal $\SL_2$-homomorphisms for $X$, each defined over $K$, then
  by Proposition \ref{prop:assoc-cochar}, the $K$-cocharacters $\Phi =
  \phi_{\mid \T}$ and $\Psi = \psi_{\mid \T}$ associated with $X$ are
  conjugate by a unique element of $U(K)$. Thus we may replace $\psi$
  by a $U(K)$-conjugate and suppose that $\phi_{\mid \T} = \psi_{\mid
    \T}$. Proposition \ref{prop:optimal-equality} then shows that
  $\phi = \psi$ and the proof is complete.
\end{proof}

\begin{rem}
  In the case of a \emph{finite} ground field $K$, Seitz
  \cite{seitz}*{Prop. 9.1} obtained existence and conjugacy over $K$
  for \emph{good $A_1$ subgroups} (see \S \ref{sub:good-a1s} below
  for their definition).
\end{rem}

\subsection{Complete reducibility of optimal homomorphisms}
\label{sub:gcr}

Let $G$ be any reductive group. Generalizing the notion of a
completely reducible representation of a group, J-P. Serre has
introduced the following definition. A subgroup $H \le G$ is said to
be $G$-completely reducible (for short: \gcr) if for every parabolic subgroup
$P$ of $G$ containing $H$ there is a Levi subgroup of $P$ which also
contains $H$. See \cite{serre-sem-bourb} for more on this notion.

We are going to prove that the image of an optimal homomorphism is
\gcr. We establish some technical lemmas needed in the proof.  First,
we show that a suitable generalization of Proposition
\ref{prop:exp-centralizer} is valid.

\begin{lem}
  \label{lem:exp-centralizer-non-Richardson}
  Let $\Psi \in X_*(G)$ and suppose that $P=P(\Psi)$ is a
  distinguished parabolic subgroup with unipotent radical $U=R_uP$.
  Suppose that the nilpotence class of $U$ is $<p$, and let
  \begin{equation*}
    \e:\Lie(U) \to U
  \end{equation*}
  be the isomorphism of Proposition \ref{prop:unipotent-exp-iso}.  If
  $X_0 \in \glie(\Psi;n)$ for some $n \ge 1$, 
   then $X_0 \in \Lie(U)$ and $C^o_G(X_0) = C^o_G(\e(X_0))$.
\end{lem}

\begin{proof}  
  Let $N(X_0) = \{g \in G \mid \Ad(g)X_0 \in k X_0\} \le G$.
  By assumption, the torus $\Psi(\G_m)$ is contained in $N(X_0)$; in
  particular, this torus normalizes $C_G(X_0)$. We may choose a
  maximal torus $S$ of $C_G(X_0)$ centralized by $\Psi(\G_m)$; thus
  $S' = S \cdot \Psi(\G_m)$ is a maximal torus of $N(X_0)$. According
  to \cite{mcninch-rat}*{Lemma 25}, there is a cocharacter $\Lambda
  \in X_*(S')$ which is associated to $X_0$. Let $T$ be a maximal torus
  of $G$ containing $S'$; thus $T$ lies in the centralizer of
  $\Lambda(\G_m)$, of $S$, and of $\Psi(\G_m)$.
 
  Since a Richardson orbit representative $X$ for the dense $P$-orbit
  on $U$ satisfies $X^{[p]}=0$, we have also $X_0^{[p]}=0$.  Now
  consider the Levi subgroup $L=C_G(S)$; the nilpotent element $X_0$
  is distinguished in $\Lie(L)$.  Let $Q=P_L(\Lambda)$, and let
  $V=R_uQ$ be the unipotent radical of $Q$.  Proposition
  \ref{prop:unipotent-exp-iso} gives a unique isomorphism
  \begin{equation*}
    \e':\Lie(V) \to V,
  \end{equation*}
  and we know from Proposition \ref{prop:exp-centralizer} that
  $C_G^o(X_0)=C_G^o(\e'(X_0))$. Thus our lemma will follow if we show
  that $\e(X_0) = \e'(X_0)$.
  
  Notice that $T$ is contained in the Levi factors $Z_G(\Psi)$ of $P$
  and $Z_L(\Lambda)$ of $Q$, so that $T$ normalizes the connected
  unipotent subgroup $W=(U \cap V)^o$ of $G$.
  Since the nilpotence class of $W$ is $<p$,
  \cite{seitz}*{Proposition 5.2} gives a unique isomorphism of algebraic groups
  \begin{equation*}
    \e'':\Lie(W) \to W
  \end{equation*}
  whose tangent map is the identity and which is compatible with the
  action of the connected solvable group $T\cdot W$ by conjugation.
  On the other hand, the tangent maps of the restrictions $\e_{\mid
    \Lie(W)}$ and $\e'_{\mid \Lie(W)}$ are the identity, and these maps
  are compatible with the action of $T \cdot W$;  we thus have
  \begin{equation*}
    \e_{\mid \Lie(W)} = \e'' = \e'_{\mid \Lie(W)}.
  \end{equation*}
  This implies that $\e(X_0)=\e'(X_0)$ as desired, and the proof is complete.
\end{proof}

We now show that a suitable deformation of an optimal homomorphism remains
optimal.
\begin{lem}
  \label{lem:optimal-lim}
  Let $\phi:\SL_2 \to G$ be an optimal $\SL_2$-homomorphism, and
  suppose that $\phi$ takes its values in the parabolic subgroup
  $P$. 
  \begin{enumerate}
  \item There is  a cocharacter $\gamma \in X_*(P)$ such that
    $\gamma(\G_m)$ centralizes $\phi(\T)$ and such that $P=P(\gamma)$.
  \item Denoting by $L=Z(\gamma)$ the Levi factor of $P$ determined by
    $\gamma$, write $\widehat\phi:\SL_2 \to L$ for the
    homomorphism 
    \begin{equation*}
      x \mapsto \lim_{t \to 0}
      \gamma(t)\phi(x)\gamma(t^{-1})
    \end{equation*}
    of Lemma \ref{lem:limit-homomorphism}. Then $\widehat \phi$ is an
    optimal $\SL_2$-homomorphism as well.
  \end{enumerate}
\end{lem}

\begin{proof}
  Since $\phi(\T)$ lies in some maximal torus of $P$, (1) follows from
  Lemma \ref{lem:parabolic-cocharacter}.
  
  Let us prove (2).  Let $X = d\phi(X_1)$ as usual, and write $\Psi$
  for the cocharacter $\phi_{\mid \T}$; it is associated with $X$.
  Denoting by $C_\Psi$ the corresponding Levi factor of the
  centralizer of $X$, we may choose a maximal torus $S \le C_\Psi$ and
  Proposition \ref{prop:optimal-sl2-Levi} implies that $\phi$ takes
  its values in the Levi subgroup $C_G(S)$.  We may evidently
  replace $G$ by $L$ and so assume that $X$ is distinguished.
  
  Now let $X = X_0 + X'$, $Y = Y_0 + Y'$ with $X_0,Y_0 \in \Lie(L) =
  \glie(\gamma;0)$ and with $X',Y' \in \Lie(R_uP)$.  Lemma
  \ref{lem:limit-homomorphism} shows that $d\widehat \phi(X_1) = X_0$
  and $d\widehat \phi(Y_1) = Y_0$.
  
  To shows that $\widehat \phi$ is optimal for $X_0$, it is enough to
  show that $\widehat \phi$ takes values in some Levi subgroup $M$ of
  $L$ such that $X_0 \in \Lie(M)$ is distinguished. Indeed, since
  $\SL_2$ is its own derived group, this will imply that $\Psi =
  \phi_{\mid \T}$ takes its values in $(M,M)$, so that $\Psi$ is
  indeed associated with $X_0$.
  
  Note that the torus $\Psi(\G_m)$ normalizes $C_L(X_0)$. Since
  $\Psi(\G_m)$ lies in a maximal torus of the semidirect product of
  $C_L(X_0)$ and $\Psi(\G_m)$, it is clear that there is a maximal
  torus $S$ of $C_L(X_0)$ centralized by $\Psi(\G_m)$. Taking
  $M=C_L(S)$, we claim that $\phi$ takes its values in $M$.
  
  Notice that 
  \begin{equation*}
    \widehat \phi(x_1(t)) = \lim_{s \to 0} \gamma(s)\e(tX)\gamma(s^{-1})= 
    \lim_{s \to 0} \e(t(X_0+\Ad(\gamma(s))X')) = \e(tX_0) 
  \end{equation*}
  for each $t \in k$, 
  Similarly, $\widehat \phi(y_1(t)) = \e(tY_0)$ for each $t \in k$.
  
  Since $S$ is contained in the centralizer of $X$, it is contained in
  the instability parabolic $P_X$ for $X$ Proposition
  \ref{prop:assoc-cochar}. Thus $\e$ is $S$-equivariant.  Since
  $\SL_2$ is generated by $\X$ and $\X^-$, this equivariance shows
  that we are done if $S$ centralizes both $X_0$ and $Y_0$ -- of
  course, $S$ centralizes $X_0$ by assumption.
  
  Write $H = d\Psi(1)$; since $\Psi$ and $\gamma$ commute, $\widehat
  \phi_{\mid \T} = \Psi$.  Now, $\ad(X_0)Y_0=[X_0,Y_0]=H$.  As in the
  proof of Proposition \ref{prop:optimal-sl2-Levi}, we write
  $Y_0=\sum_{\lambda \in X^*(S)} Y_{0,\lambda}$ as a sum of weight
  vectors for the torus $S$. Since $\Psi(\G_m)$ commutes with $S$, $H$
  is centralized by $S$, and so we have $[X_0,Y_{0,\lambda}]=0$ when
  $\lambda \not = 0$; we want to conclude that $Y_{0,\lambda}=0$. We
  do not know that $\Psi$ is associated with $X_0$, so we can not
  simply invoke Proposition \ref{prop:assoc-cochar}.  However, since
  $Y_{0,\lambda} \in \glie(\Psi;-2)$, the general theory of
  $\SL_2$-representations shows: if $Y_{0,\lambda}\not =0$, then
  $\widehat \rho(x_1(t)) = \e(tX_0)$ acts non-trivially on
  $Y_{0,\lambda}$ for some $t \in k^\times$.  On the other hand,
  according to Lemma \ref{lem:exp-centralizer-non-Richardson} we have
  $C_L^o(X_0) = C_L^o(\e(tX_0))$, so that $Y_{0,\lambda} \in
  \clie_{\Lie(L)}(X_0)= \clie_{\Lie(L)}(\e(tX_0))$.  Thus indeed
  $Y_{0,\lambda}=0$ for each non-zero $\lambda$, as required.  Thus
  $Y_0=Y_{0,0}$ so that $S$ centralizes $Y_0$; the proof is now
  complete.
\end{proof}

\begin{lem}
  \label{lem:recognize-assoc}
  Let $X \in \glie$ be any nilpotent element, let $\psi
  \in X_*(G)$ a cocharacter associated with $X$, and let
  $L=C_G(\psi(\G_m))$ be the Levi factor in the instability parabolic
  determined by $\psi$.
  \begin{enumerate}
  \item The $L$  orbit $\mathcal{V} = \Ad(L)X$ is
    a Zariski open subset of $\glie(\psi;2)$.
  \item Let $Y \in \glie$ be nilpotent. Then $\psi$ is a cocharacter
    associated with $Y$ if and only if $Y \in \mathcal{V}$.
  \end{enumerate}
\end{lem}

\begin{proof}
  To prove (1), note that the orbit map 
  \begin{equation*}
    y \mapsto \Ad(y)X:L \to
  \glie(\psi;2)
  \end{equation*}
  has differential $\ad(X):\Lie(L)=\glie(\psi;0) \to \glie(\psi;2)$;
  if we know that the differential is surjective, then the orbit map
  is dominant and separable and (1) follows.  To see the surjectivity,
  we argue as follows. Recall from Proposition
  \ref{prop:assoc-cochar} that $\clie_\glie(X)$ is
  contained in $\sum_{i \ge 0} \glie(\psi;i)$; in particular,
  $\glie(\psi;-2) \cap \clie_\glie(X)=0$. According to
  \cite{jantzen-nil}*{Lemma 5.7} this last observation implies (in
  fact: is equivalent to) the statement $[\glie(\psi;0),X] =
  \glie(\psi;2)$; this proves the required surjectivity (note that
  \cite{jantzen-nil}*{5.7} is applicable since the Lie algebra of a
  strongly standard reductive group has on it a nondegenerate,
  invariant, symmetric, bilinear form -- cf.  Proposition
  \ref{prop:strongly-standard}).
  
  For (2) note first that $\psi$ is evidently associated to any $Y \in
  \mathcal{V}$. Conversely, if $\psi$ is associated to $Y$, then $Y
  \in \glie(\psi;2)$, and (1) shows that $\Ad(L)Y$ is also open and
  dense in $\glie(\psi;2)$. Thus $\Ad(L)X \cap \Ad(L)Y \not =
  \emptyset$, so that $Y \in \Ad(L)X = \mathcal{V}$.
\end{proof}

\begin{theorem}
  \label{theorem:optimal-gcr}
  Let $G$ be strongly standard, and let $\phi:\SL_2 \to G$ be an
  optimal $\SL_2$ homomorphism. Then the image of $\phi$ is \gcr.
\end{theorem}

\begin{proof}
  Let $X = d\phi(X_1)$ as usual, and write $\Psi$ for the cocharacter
  $\phi_{\mid \T}$; it is associated with $X$.  Denoting by $C_\Psi$
  the corresponding Levi factor of the centralizer of $X$, we may
  choose a maximal torus $S \le C_\Psi$ and Proposition
  \ref{prop:optimal-sl2-Levi} implies that $\phi$ takes its values in
  the Levi subgroup $L=C_G(S)$. Applying \cite{serre-sem-bourb}*{Prop.
    3.2}, one knows that $\phi(\SL_2)$ is \gcr~ if and only if it is
  $L$-\text{cr}.
%  %%
%  \footnote{The following result is due to J-P. Serre: If $L < G$ is a
%    Levi subgroup, a subgroup $H<L$ is $L$-$\operatorname{cr}$ if and
%    only if $H$ is \gcr. A proof was given by Serre in a 1998 course
%    at Princeton University.}.
%  %%
  We replace $G$ by $L$, and thus suppose that $X$ is
  \emph{distinguished.}
  
  Let $P$ be a parabolic subgroup of $G$ and suppose that the image of
  $\phi$ lies in $P$.  We claim that since $X$ is distinguished, we
  must have $P=G$; this will prove the theorem.
  
  To prove our claim, first notice that by Lemma
  \ref{lem:optimal-lim}(1) we may choose $\gamma \in X_*(P)$ with
  $P=P(\gamma)$ and such that $\gamma(\G_m)$ commutes with
  $\Psi(\G_m)$.
  
  Let us write $X = \sum_{i \ge 0} X_i$ with $X_i \in
  \glie(\gamma;i)$.  Consider the homomorphism $\widehat \phi:\SL_2
  \to Z(\gamma)$ constructed in Lemma \ref{lem:optimal-lim}; according
  to (2) of that lemma, $\widehat \phi$ is optimal for $X_0$, so that
  the cocharacter $\Psi$ is associated to $X_0$ as well as to $X$.
  
  We now claim that $X$ and $X_0$ are conjugate. This will show that
  $X_0$ is distinguished in $G$, hence that $G=Z(\gamma)$ so that also
  $P=G$ as desired. Let $L=C_G(\Psi(\G_m))$. Then Lemma
  \ref{lem:recognize-assoc} implies that $X_0$ is contained in the
  orbit $\mathcal{V}=\Ad(L)X \subset \glie(\Psi;2)$, proving our
  claim.
\end{proof}

\subsection{Comparison with good homomorphisms}
\label{sub:good-a1s}

According to Seitz \cite{seitz}, an $\SL_2$ homomorphism $\phi:\SL_2 \to G$
is called \emph{good} (or \emph{restricted}) provided that the weights of
a maximal torus of $\SL_2$ on $\Lie(G)$ are all $\le 2p-2$. 

\begin{prop}
  Let $\phi:\SL_2 \to G$ be a homomorphism, where $G$ is a strongly
  standard reductive group. Then $\phi$ is good if and only if it is
  optimal for $X = d\phi(X_1)$.  In particular, all good
  $\SL_2$-homomorphisms whose image contains the unipotent element $v$
  are conjugate by $C_G^o(v)$.
\end{prop}

\begin{proof}
  That an optimal homomorphism is good follows from Proposition
  \ref{prop:p-nilpotent-cocharacter}.  Choose a Springer isomorphism
  $\Lambda:\UU \to \NN$. If $u$ is a unipotent element of order $p$,
  choose a Levi subgroup $L$ in which $u$ is distinguished; this just
  means that $X=\Lambda(u) \in \glie$ is distinguished. It follows
  from Proposition \ref{prop:order-result} that $X^{[p]}=0$. Choose an
  optimal homomorphism $\phi'$ for $X$; we know that $\phi'$ takes
  values in $L$ (Proposition \ref{prop:optimal-sl2-Levi}), and if $v =
  \phi'(x(1))$, it is clear from Proposition
  \ref{prop:sl2-exponential-map} that $v$ and $u$ are Richardson
  elements in the same parabolic subgroup of $L$; thus $v$ and $u$ are
  conjugate. This proves that $u$ is in the image of some optimal
  homomorphism $\phi$.

  To prove that good homomorphisms are optimal, we use a result of
  Seitz. Since $\phi$ is optimal, we just observed that it is good,
  and Seitz proved \cite{seitz}*{Theorem 1.1} that any good
  homomorphism with $u$ in its image is conjugate by $C_G(u)$ to
  $\phi$. Thus, any good homomorphism is indeed optimal.
\end{proof}

\section{Rational elements of a nilpotent orbit defined over a ground field}
\label{sec:kottwitz}

In this section, we extend a result first obtained by R. Kottwitz
\cite{kottwitz} in the case where $K$ has characteristic 0.  We give
here a proof which is also valid in positive characteristic (under
some assumptions on $G$). For the most part, we follow the original
argument of Kottwitz.
\begin{theorem}
  Let $K$ be any field, and let $G$ be a strongly standard
  connected reductive $K$-group which is $K$-quasisplit.  If the nilpotent orbit
  $\mathcal{O} \subset \NN$ is defined over $K$, then $\mathcal{O}$
  has a $K$-rational point.
\end{theorem}

\begin{proof}
  If $K$ is a finite field, the theorem is a consequence of the
  Lang-Steinberg theorem; cf. \cite{steinberg-endomorphisms}*{\S10}
  and \cite{steinberg-regular}.  Suppose now $K$ to be infinite.

  We fix a Borel subgroup $B$ of $G$ which is defined over $K$, and a
  maximal torus $T \subset B$ which is also over $K$.  The roots of
  $G$ in $X^*(T)$ which appear in the Lie algebra of the unipotent
  radical of $B$ are declared positive, and we will write $\overline{C} \subset
  X_*(T)$ for the positive Weyl chamber determined by $B$:
  \begin{equation*}
    \overline{C} = \{\mu \mid \langle \alpha,\mu \rangle \ge 0 
    \ \text{for all positive roots $\alpha$ of $G$ in $X^*(T)$}\}.
  \end{equation*}
  If $W = N_G(T)/T$ denotes the Weyl group of $T$, then each $\mu \in
  X_*(T)$ is $W$-conjugate to a unique point in $\overline{C}$.  We
  also write $\Gamma = \Gal(K_\sep/K)$ for the absolute Galois group
  of the field $K$.

  The $K$-variety $\mathcal{O}$ has a point $X'$ rational over the
  separable closure $K_\sep$ of $K$ in $k$ (e.g. by
  \cite{springer-LAG}*{11.2.7}). According to Proposition
  \ref{prop:assoc-cochar}, there is a cocharacter $\Psi'$ associated
  with $X'$ and defined over $K_\sep$. Let $T'$ be a maximal torus of
  $G$ defined over $K_\sep$ which contains the image of $\Psi'$.
  
  For $\gamma \in \Gamma$, the cocharacter ${\Psi'}^\gamma$ is
  associated with the nilpotent ${X'}^\gamma$. Since $\mathcal{O}$ is
  defined over $K$, ${X'}^\gamma$ and $X'$ are conjugate. Hence
  $\Psi'$ and ${\Psi'}^\gamma$ are conjugate by another application of
  Proposition \ref{prop:assoc-cochar}.
  
  According to \cite{springer-LAG}*{Prop. 13.3.1 and 11.2.7} we may
  find $g \in G(K_\sep)$ such that $gT'g^{-1} = T$; the same reference
  shows that any element $w$ of the Weyl group of $T$ may be
  represented by an element $\dot w \in N_G(T)$ rational over
  $K_\sep$.  We have that $\Psi = \Int(g) \circ \Psi' \in X_*(T)$ is
  defined over $K_\sep$.  Replacing $\Psi$ by $\Int(\dot w) \circ
  \Psi$ for a suitable $w$ in the Weyl group of $T$, we may suppose
  that $\Psi \in \overline{C} \subset X_*(T)$ and is defined over
  $K_\sep$. Of course, $\Psi$ is associated with the nilpotent element
  $X = \Ad(\dot w g)X'$.
  
  Since $B$ and $T$ are $\Gamma$-stable, $\gamma$ permutes the
  positive roots in $X^*(T)$. Thus, $\gamma$ leaves $\overline{C}$
  invariant; in particular, $\Psi^\gamma \in \overline{C}$.  We know
  $\Psi$ and $\Psi^\gamma$ to be conjugate in $G$. Since $T$ is a
  maximal torus of the centralizer of both $\Psi(\G_m)$ and of
  $\Psi^\gamma(\G_m)$, we may suppose that $\Psi^\gamma = \Int(\dot
  w)\Psi$ for some $w$ in the Weyl group of $T$. But $\overline{C}$ is
  a fundamental domain for the $W$-action on $X_*(T)$, so we see that
  $\Psi = \Psi^\gamma$. Since $\Psi$ is defined over $K_\sep$ and is
  $\Gamma$-stable, $\Psi$ is defined over $K$
  \cite{springer-LAG}*{11.2.9}.
  
  This shows in particular that the subspace $\glie(\Psi;2)$ is
  defined over $K$. According to Lemma \ref{lem:recognize-assoc},
  there is a Zariski open subset of $\glie(\Psi;2)$ consisting of
  elements in $\mathcal{O}$. Since $K$ is infinite, the $K$-rational
  points of $\glie(\Psi;2)$ are Zariski dense in $\glie(\Psi;2)$. Hence
  there is a $K$-rational point in $\mathcal{O}$ and the proof is
  complete.
\end{proof}

\begin{cor}
  Let $G$ be a strongly standard reductive $K$-group which is
  $K$-quasisplit.  There is a regular nilpotent element $X \in
  \glie(K)$.  In particular, there is an optimal homomorphism
  $\phi:\SL_2 \to G$ defined over $K$ with $d\phi(X_1)=X$.
\end{cor}

\begin{proof}
  Since $G$ is split over a separable closure $K_\sep$ of $K$, there
  is a $K_\sep$ rational regular nilpotent element. Thus the regular
  nilpotent orbit is defined over $K_\sep$. Since this orbit is
  clearly stable under $\Gal(K_\sep/K)$, it is defined over $K$. So
  the theorem shows that there is a $K$-rational regular nilpotent
  element $X$. The final assertion follows from Theorem
  \ref{theorem:optimal-over-ground-field}.
\end{proof}

\begin{rem}
  With $G$ as in the theorem, there is a Springer isomorphism
  $\Lambda:\UU \to \NN$ defined over $K$. Thus a unipotent
  conjugacy class defined over $K$ has a $K$-rational point.
\end{rem}

\section{Appendix: Springer Isomorphisms (Jean-Pierre Serre, June 1999)}

\newcommand{\Gu}{\ensuremath{G^{\operatorname{u}}}}
\newcommand{\Gur}{\ensuremath{G^{\operatorname{ur}}}}
\newcommand{\Cur}{\ensuremath{C(u)^{\operatorname{r}}}}
\newcommand{\glien}{\ensuremath{\glie^{\operatorname{n}}}}
\newcommand{\glienr}{\ensuremath{\glie^{\operatorname{n}}}}

%% re-init the ``theorem styles''; Serre's notes didn't have numbered theorems.
\theoremstyle{plain}
\newtheorem*{proposition}{Proposition}
\newtheorem*{corollary}{Corollary}
\newtheorem*{lemma}{Lemma}

\theoremstyle{remark}
\newtheorem*{note}{Note}

Let $G$ be a simple algebraic group in char. $p$, which I assume to be
``good'' for $G$. I also assume the ground field $k$ to be
algebraically closed. Call $\Gu$ the variety of unipotent elements of
$G$ and $\glie^n$ the subvariety of $\glie = \Lie(G)$ made up of the
nilpotent elements.

Springer has shown that there exist algebraic morphisms 
\begin{equation*}
  f:\Gu \to \glien
\end{equation*}
with the following properties:
\begin{itemize}
\item[a)] $f$ is compatible with the action of $G$ by conjugation on both sides. 
\item[b)] $f$ is bijective.
\end{itemize}
In fact, it was later shown that these properties imply (at least when
$p$ is ``very good'', which is always the case if $G$ is not of type
$A$):
\begin{itemize}
\item[b$'$)] $f$ is an isomorphism of algebraic varieties.
\end{itemize}

Despite the fact that there are \emph{many such} $f$'s (they make up
an algebraic variety of dimension $\ell$, where $\ell$ is the rank of
$G$), one often finds in the literature the expression ``the Springer
isomorphism'' used --and abused --, especially to conclude that the
$G$-classes of unipotent elements of $G$ and nilpotent elements of
$\glie$ are in a natural correspondence, namely ``the'' Springer
correspondence.  

It might be good for the reader to consider the case of $G=\SL_n$ (or
rather $\PGL_n$, if one wants an adjoint group). In that case a
Springer isomorphism is of the form
\begin{equation*}
  1 + e \mapsto a_1e + \cdots + a_{n-1}e^{n-1},
\end{equation*}
where $e^n = 0$ (so that $u=1+e$ is unipotent), and the $a_i$ are
elements of $k$ with $a_1 \not = 0$. Every such family $\vec
a=(a_1,\dots,a_{n-1})$ defines a unique Springer isomorphism $f_{\vec
  a}$, and one gets in this way every Springer isomorphism, once and
only once.  This example also shows that the Springer isomorphisms can
be quite different: e.g., for some one may have $f(u^m) = m.f(u)$ for
all $u$ and all $m\in \Z$ ( such an $f$ exists if and only if $p \ge
n)$, and for some one does not even have $f(u^{-1}) = -f(u)$!

In what follows, I want to repair this unfortunate mix-up by showing
that all the different Springer isomorphisms give \emph{the same}
bijection between the $G$-classes of $\Gu$ and the $G$-classes of
$\glien$, so that one can indeed speak (in that case) of \emph{the}
Springer bijection.

I have to recall first how the Springer isomorphisms are defined. Call
$\Gur$ the set of regular unipotent elements of $G$; it is an open
dense set in $\Gu$; same definition for $\glienr$ in $\glie =
\Lie(G)$.  Choose an element $u$ in $\Gur$ and let $C(u)$ be its
centralizer. It is known that $C(u)$ is smooth, connected, unipotent,
commutative, of dimension $\ell$ (= rank $G$). Let $\clie(u) = \Lie
C(u)$ be its Lie algebra.  Choose an element $X$ of $\clie(u)$ which
is regular. Then its centralizer is $C(u)$, and the Springer
construction shows that there is a \emph{unique} Springer isomorphism
$f = f_{u,X}$ which has the property that $f(u) = X$. Let us fix $X$;
then it is clear that every Springer isomorphism is equal to $f_{v,X}$
for some $v \in \Cur$, where $\Cur = C(u) \cap \Gur$; moreover, $v$ is
uniquely defined by $f$. Hence we have a \emph{one-to-one
  parametrization of the Springer isomorphisms by the elements $v$ of
  $\Cur$.}

The next step consists in showing that this parametrization is ``algebraic''.
The precise meaning of this is the following:

\begin{proposition}
  There exists an algebraic morphism $F:\Cur \times \Gu \to \glien$
  such that $F(v,z) = f_{v,X}(z)$ for every $v \in \Cur$ and $z \in
  \Gu$.
\end{proposition}

\begin{proof}
  Call $N_u$ the normalizer of $C(u)$ in $G$. Since all regular
  unipotents are conjugate, $N_u$ acts transitively on $\Cur$, so that
  one can identify the algebraic variety $\Cur$ with the coset space
  $N_u/C(u)$.  Similarly, one may identify $\Gur$ with $G/C(u)$. Let
  us now define an algebraic map
\begin{equation*}
  F':N_u \times G \to \glien
\end{equation*}
by the formula 
\begin{equation*}
  F'(n,z) = \Ad(zn^{-1}).X 
\end{equation*}
(i.e. the image of $X \in \glie$ by the inner automorphism defined by
$zn^{-1}$).  It is clear that $F'(n,z)$ depends on $n$ only mod.
$C(u)$, and that it depends on $z$ also mod $C(u)$. Hence $F'$ factors
out and gives a map of $N_u/C(u) \times G/C(u)$ into $\glien$. If we
identify $N_u/C(u)$ with $\Cur$ and $G/C(u)$ with $\Gur$, we thus get
a map
\begin{equation*}
  F_0:\Cur \times \Gur \to \glien.
\end{equation*}
It is well-known that $\Gu$ is a normal variety and that $\Gu
\setminus \Gur$ has codimension $> 1$ in $\Gu$. Hence the same is true
for $\Cur \times \Gur$ in $\Cur \times \Gu$. Since $\glien$ is an
affine variety, the map $F_0$ extends uniquely to an algebraic map
$F:\Cur \times \Gu \to \glien$.  One checks immediately that for every
fixed $v \in \Cur$, the map $z \mapsto F(v,z)$ has the following
properties: a) it commutes with the action of $G$; b) it maps $v$ to
$X$.  (Property a) is checked on $\Gur$ first; by continuity, it is
valid everywhere.) This shows that $F$ is the map we wanted.
\end{proof}

\begin{corollary}
  The bijection 
  \begin{equation*}
    \text{$G$-classes of $\Gu$} \to \text{$G$-classes of $\glien$}
  \end{equation*}
  given by a Springer isomorphism $f$ is independent of the choice of
  $f$.
\end{corollary}

This is easy. One uses the following elementary lemma:
\begin{lemma}
  Let $Y$, $Z$ be two $G$-spaces. Assume $G$ has finitely many orbits
  in each.  Let $T$ be a connected space, and $F:T \times Y \to Z$ a
  morphism such that, for every $t \in T$, the map $y \mapsto F(t,y)$
  is a $G$-isomorphism of $Y$ on $Z$.
  
  Then, for every $y \in Y$, the points $F(t,y)$, $t \in T$, belong to
  the same $G$-orbit.
\end{lemma}

Proof by induction on $\dim Y = \dim Z$. The statement is clear in
dimension zero, because of the connexity of $T$. If $\dim Y > 0$,
there are finitely many open orbits in $Y$ (resp. $Z$); call $Y_0$ and
$Z_0$ their union. It is clear that, for every $t$, the isomorphism
$F_t:y \mapsto F(t,y)$ maps $Y_0$ into $Z_0$. Moreover, the connexity
of $T$ implies that the $F_t$'s map a given connected component of
$Y_0$ into the same connected component of $Z_0$. And the induction
hypothesis applies to $Y \setminus Y_0$ and $Z \setminus Z_0$.

The corollary follows from the lemma, applied with $T = \Cur$, $Y=
\Gu$ and $Z = \glien$.

\begin{note}
  The structure of $N_u/C(u)$ seems interesting. If I am not mistaken,
  it is the semi-direct product of $\G_m$ by a unipotent connected
  group $V$ of dimension $\ell - 1$; moreover, the action of $\G_m$ on
  $\Lie V$ has weights equal to $k_2-1,k_3-1,\dots,k_\ell -1$, where
  the $k_i$'s are the exponents of the Weyl group.
  
  Another interesting (and related) question is the behaviour of a
  Springer isomorphism $f$ when one restricts $f$ to $C(u)$. The
  tangent map to $f$ is an endomorphism of $c(u) = \Lie C(u)$. Is it
  always a non-zero multiple of the identity?
\end{note}

\vspace{5mm}
\quad J-P. Serre \quad\quad June 1999

\newcommand{\myresetbiblist}[1]{%
  \settowidth{\labelwidth}{\def\thebib{#1}\BibLabel}%
  \setlength\labelsep{1mm}
  \setlength\leftmargin\labelwidth
  \addtolength\leftmargin\labelsep
}

\newcommand\mylabel[1]{#1\hfil}

%\bibliography{MathBib,books,Preprints}

\end{document}